\documentclass[review,hidelinks,onefignum,onetabnum]{siamart220329}

\usepackage{amsfonts}
\usepackage{graphicx}
\usepackage{epstopdf}
\usepackage{algpseudocode}
\usepackage{subcaption}

\ifpdf
  \DeclareGraphicsExtensions{.eps,.pdf,.png,.jpg}
\else
  \DeclareGraphicsExtensions{.eps}
\fi

\newsiamremark{remark}{Remark}
\newsiamremark{hypothesis}{Hypothesis}
\crefname{hypothesis}{Hypothesis}{Hypotheses}
\newsiamthm{claim}{Claim}

\headers{Globally Convergent Homotopies for Optimal Control}{W. Esterhuizen, K. Flaßkamp, M. Hoffmann, K.Worthmann}

\title{Globally Convergent Homotopies for Discrete-Time Optimal Control\thanks{Submitted to the editors DATE.
\funding{This work is funded by the Deutsche Forschungsgemeinschaft (DFG), project number 461953135}}}

\author{Willem Esterhuizen\footnotemark[2]\  \footnotemark[3]
\and Kathrin Flaßkamp\thanks{Systems Modeling and Simulation, Systems Engineering, Universität des Saarlandes, Germany
  (\email{kathrin.flasskamp@uni-saarland.de}, \email{matthias.hoffmann@uni-saarland.de}).}
\and Matthias Hoffmann\footnotemark[2]
\and Karl Worthmann\thanks{Optimization-based Control Group, Institute of Mathematics, Technische Universität Ilmenau, Germany
	(\email{willem-daniel.esterhuizen@tu-ilmenau.de}, \email{karl.worthmann@tu-ilmenau.de}).}}

\usepackage{amsopn}

\usepackage{graphicx}      
\usepackage{amsmath} 
\usepackage{amssymb}  
\usepackage{color}
\usepackage{graphicx}
\usepackage{booktabs} 
\usepackage{cases}
\usepackage{tikz}
\usepackage{pgfplots}
\pgfplotsset{compat=newest}
\usepackage{siunitx}

\usetikzlibrary{arrows,shapes,positioning,calc,intersections,through,backgrounds,fit,decorations.markings, arrows.meta}
\usetikzlibrary{spy}
\usepackage{adjustbox}
\usepgfplotslibrary{groupplots}

\newcommand{\mbbR}{\mathbb{R}}

\newcommand{\OCP}{\mathrm{OCP}}
\newcommand{\NLP}{\mathrm{NLP}}

\newcommand{\boldu}{\mathbf{u}}

\newcommand{\rank}{\mathrm{rank}}
\newcommand{\x}{\boldsymbol x}
\newcommand{\zero}{\boldsymbol 0}

\makeatletter
\newcommand*\bigcdot{\mathpalette\bigcdot@{0.75}}
\newcommand*\bigcdot@[2]{\mathbin{\vcenter{\hbox{\scalebox{#2}{$\m@th#1\bullet$}}}}}
\makeatother

\begin{document}

\nolinenumbers

\maketitle

\begin{abstract}
    Homotopy methods are attractive due to their capability of solving difficult optimisation and optimal control problems.
    The underlying idea is to construct a homotopy, which may be considered as a continuous (zero) curve between the difficult original problem and a related, comparatively easy one. Then, the solution of the easier one is continuously perturbed along the zero curve towards the sought-after solution of the original problem. 
    We propose a methodology for the systematic construction of such zero curves for discrete-time optimal control problems drawing upon the theory of globally convergent homotopies for nonlinear programs. 
    The proposed framework ensures that for almost every initial guess at a solution there exists a suitable homotopy path
    that is, in addition, numerically convenient to track. 
    We demonstrate the results by solving optimal path planning problems for a linear system and the nonlinear nonholonomic car (Dubins' vehicle). 
\end{abstract}

\begin{keywords}
	homotopy methods, optimal control, nonlinear programming, path planning
\end{keywords}

\begin{MSCcodes}
	49K15, 
	49M99, 
	90C30, 
	93B40, 
	93C55 
\end{MSCcodes}

\section{Introduction}\label{sec:introduction}

If faced with a difficult nonlinear and, in addition, potentially non-convex discrete-time optimal control problem (OCP), one may resort to \emph{homotopy methods} for its solution.
The basic idea is to identify the key source of difficulty and, then, to introduce a \emph{homotopy parameter} $\lambda \in [0,1]$ such that $\lambda = 0$ corresponds to a significantly easier-to-solve problem while the original problem is recovered for $\lambda = 1$. 
One then continuously perturbs (in some way) the solution of the easy problem towards the solution of the difficult problem.  

A discrete-time finite-horizon OCP is a nonlinear program (NLP), where the constraints and, potentially, the objective function explicitly depend on the state and, thus, implicitly on the decision variable (the control function), via the system dynamics. 
Hence, if we are interested in questions concerning the \emph{convergence} of homotopy algorithms applied to discrete-time optimal control, then we should look to the nonlinear programming literature.
A typical approach to solve 
an NLP via homotopy is to consider the first order Karush-Kuhn-Tucker (KKT) necessary conditions, see for example \cite[Ch. 5]{boyd2004convex} \cite[Ch. 3]{bertsekas1999nonlinear}, and to express them as an equivalent system of equations involving the primal and dual variables as decision variables. 
This equivalent system can be derived in several ways, see for example \cite{mangasarian1976equivalence}, \cite[Ch. 4]{zangwil1981pathways} and \cite{papalambros1989nonlinear}.
By introducing the homotopy parameter~$\lambda$ into the NLP, this becomes a \textit{parametrized} system of equations. 
Then, one numerically tracks the so-called \emph{zero curve} (that is, the set of points for which the system vanishes) starting at~$\lambda = 0$, see 
\cite[Ch. 4]{zangwil1981pathways} and \cite{papalambros1989nonlinear} as well as 
\cite{allgower2003introduction} and \cite{rheinboldt2000numerical} for the general theory as well as appropriate numerical methods for tracking the zero curve.

The extensive work by Watson et al., see, e.g., \cite{watson1979solving, watson2001theory, watson2002probability, watson2013globally} to name but a few, elaborates a theory of homotopies for mathematical programming problems that are \emph{globally convergent with probability one}.
This means that, for almost all 
starting points (w.r.t.\ 
the Lebesgue measure), 
there exists a zero curve to the homotopy that successfully reaches $\lambda = 1$ and, thus, solves the original problem. 
This body of work utilises the fundamental theory derived in the seminal paper~\cite{chow1978finding} by Chow, Mallet-Paret and Yorke, in which convergent homotopies for finding zeros to maps were considered. 
Nonconvex NLPs as considered in~\cite{watson2001theory} are of particular interest to our setting of discrete-time optimal control.
However, it can be quite challenging to verify some of the key assumptions of 
the 
main result (\cite[Thm. 7.1]{watson2001theory}) when considering the OCP by itself because the conditions 
are stated at the level of the NLP. 

In this paper we make the following contributions. 
First, for NLPs, we show that a new sufficient condition implies that a difficult-to-verify assumption from \cite[Thm. 7.1]{watson2001theory} holds (stating that the homotopy's zero curve cannot loop back to $\lambda = 0$ as you traverse along it). 
We then show that under certain assumptions a particular homotopy (which approaches the KKT necessary conditions as $\lambda$ approaches 1, thereby finding a local solution to an NLP) is rendered globally convergent. 
Next, we explain how these assumptions on the NLP may be ``translated'' to assumptions on the discrete-time OCP from where the NLP originates. 
Our main result states that if all functions in the OCP are $C^3$; the original problem, where $\lambda = 1$, has a feasible solution; the control is bounded; $\lambda$ is introduced such that the nonconvex constraints are trivially satisfied with $\lambda = 0$, and if the feasible space shrinks as $\lambda$ increases, then the specified homotopy is rendered globally convergent. 

The main motivation for this research comes from path planning problems for dynamical systems through obstacle fields.
Approaches have been proposed where a homotopy parameter perturbs the obstacles, see for example \cite{kontny2016fast, bergman2018combining, mazzia2022minimum, ferbach1998method}, and our recent paper~\cite{hoffmann2023path}. 
However, the authors do not present convergence guarantees for their algorithms.

Note that the current paper’s \emph{direct} setting differs from the \emph{indirect} one considered in most papers concerned with homotopy methods in optimal control. There the focus is usually on finding zeros to the so-called \emph{shooting function} associated with a boundary value problem obtained from a study via Pontryagin's maximum principle, see the papers \cite{trelat2017geometric, caillau2012differential, caillau2012minimumtime, gergaud2006homotopy} to name a few. 
Though some work on homotopies for problems with state constraints exists for the indirect setting, 
see \cite{hermant2010homotopy,hermant2011optimal}, it would not be suitable for the mentioned application of path planning due to the presence of an extremely high number of constraints. 
(In \cite{hoffmann2023path} the real word example had over one thousand obstacles).
Optimal control problems with state constraints are notoriously difficult to analyse from an indirect perspective and the inclusion of a perturbing homotopy parameter only complicates the analysis further. 
In particular, so-called \emph{active arcs} may appear and disappear as the parameter is varied. 

\subsection{Outline}

We conclude this introduction with brief summaries of other work related to homotopy methods for solving NLPs, in Subsection~\ref{subsec:lit_review}. 
In Section~\ref{sec:homo_for_sys_eqns} we then present a summary of the well-known theory of globally convergent homotopies as it applies to solving systems of equations.  Section~\ref{sec:homoh_for_NLP} shows how the content of Section~\ref{sec:homo_for_sys_eqns} can be applied to solve nonconvex NLPs. 
The first contribution of the paper appears in Subsection~\ref{subsec:new_conditions}, where we show that with a new assumption (saying that the nonconvex constraints are trivially satisfied when $\lambda = 0$) a specific homotopy can be rendered globally convergent. 
In Section~\ref{sec:homo_for_OCPs} we consider OCPs, and show how the conditions of the main result from Section~\ref{sec:homoh_for_NLP} can be guaranteed with assumptions on the OCP. 
Section~\ref{sec:numerics} presents some numerical examples to clarify the paper's theory and to demonstrate how the results may be used in practice.
Section~\ref{sec:conclusion} concludes the paper with some comments and suggestions for future research. 

\subsection{Related work on homotopies for nonlinear programs}\label{subsec:lit_review}

We provide a brief overview 
on homotopy methods applied to nonlinear programming ---~ with a particular focus on 
convergence. 
In~\cite{poore1988expanded}, the authors propose 
an unconstrained optimisation problem where the cost function, consisting of the original cost and penalty functions involving the constraint functions, serves as a homotopy. 
Then by driving the homotopy parameter (involved in the weighting of the penalty functions) to zero, a solution is found to the original problem.
The main result states that, for certain classes of penalty functions, it is guaranteed that there exists a zero path of the homotopy locally about the desired final point. 
This work may also be interesting to our setting of convergent homotopies for optimal control, though it has been reported to be tricky to use in practice due to ill-conditioning issues~\cite{watson1990survey}.  
In the papers~\cite{zhenghua1996combined,feng1998existence}, the authors consider a homotopy formed via the KKT conditions that ignores the inequality constraints by imposing a \emph{normal cone condition}, which states that the normal cone of the feasible set does not intersect the feasible set's interior. 
Then, the results from~\cite{chow1978finding} are invoked to argue that this homotopy is globally convergent.
For results of this approach applied to problems with a very large number of constraints relative to the dimension of the decision variable, see~\cite{zhou2014flattened}. Similar ideas tailored to nonlinear semi-definite programs are considered in~\cite{yang2013homotopy}.

In~\cite{bates2007numerical,rostalski2011numerical}, the authors draw upon theory for solving systems of polynomial equations via homotopy, see also \cite{sommese2005numerical}, to solve discrete-time OCPs.
Their approach consists of two steps.
The first step, which is computationally demanding, 
solves the equality constraints of the KKT conditions for a large number of initial states via a homotopy.
The second step then finds the solution to the full KKT system via a coefficient parameter homotopy. 
This second step is computationally light and may be invoked as soon as the true current state is available, making it attractive as a type of explicit model predictive control (MPC), see for example \cite{tondel2003algorithm}. 
Since the focus is on polynomial systems, there are some powerful results concerning the existence of zero curves of the considered homotopies, see \cite[Ch. 7]{sommese2005numerical}. 

In~\cite{dunlavy2005homotopy}, an algorithm for globally solving constrained optimisation problems via a homotopy on the cost function is presented. 
First, 
the homotopy parameter perturbs the cost function from one with a trivial solution towards the desired cost function. 
Then, the resulting local solution 
is used as a warm start for the next problem, while the homotopy parameter monotonically increases at each step of this iterative procedure. 
This algorithm is then extended so that at each update of the parameter a larger number of initial guesses is created by perturbing the previous local solutions. 
The algorithm then finds a number of local solutions. 
The authors show that the approach is globally convergent, in the sense of~\cite{chow1978finding,watson2001theory}, if the cost function is convex and the decision space is unconstrained.
In~\cite{best1996algorithm}, the authors consider parametric quadratic programs (QPs) and introduce the \emph{parametric quadratic programming algorithm}, which, in fact, is a \emph{parametric active set method}, see, e.g., \cite{ferreau2014qpoases}. 
In the algorithm, the monotonically increasing (homotopy) parameter 
perturbs the solution to a previous QP in the direction of the solution of the next QP (with the updated parameter). 
The primal or dual variables along with the active set of constraints are adapted at each iteration, if needed, in order to have feasibility of the updated solution candidate. 
We refer to~\cite{ferreau2014qpoases} for details of a software toolbox along with in-depth coverage of the approach and improvements. 
The approach is tailored to convex QPs, but may be used to find critical points of nonconvex QPs as well.

%


Finally, in~\cite{potschka2021sequential}, the setting is abstract Hilbert spaces, with the decision variable constrained to a convex set and subjected to equality constraints, making it particularly relevant for solving OCPs involving ordinary or partial differential equations.
The authors' \emph{sequential homotopy method} involves using homotopy to solve subproblems obtained after each Euler step of a particular \emph{gradient/anti-gradient flow}, which they introduce.
Because the state-control pair is constrained to a convex set we cannot use their current theory to solve our desired path planning problems where the obstacles are nonconvex with respect to the state.

\subsection*{Notation}

\sloppy
A function is said to be of class $C^k$, $k\in\mathbb{N}_0$, if it has a $k$-th derivative that is continuous over its domain. 
Consider a $C^1$ function $f:\mbbR^n\times\mbbR^m\rightarrow\mbbR^p$, where $f:= (f_1,f_2,\dots,f_p)^\top$, along with $x := (x_1,x_2,\dots,x_n)^\top$ and $y := (y_1,y_2,\dots,y_m)^\top$. 
Then $\nabla f(\bar x, \bar y)\in\mbbR^{p\times n}$ denotes the Jacobian of $f$ at the point $(\bar x, \bar y)\in\mbbR^{n}\times\mbbR^m$. 
Moreover, $\nabla_xf(\bar x,\bar y)\in\mbbR^{p \times n}$ denotes the Jacobian of $f$ with respect to $x$, explicitly,
\[
	\nabla_xf(\bar x,\bar y) = 
	\left(
	\begin{array}{cccc}
		\frac{\partial f_1}{\partial x_1}(\bar x, \bar y) & \frac{\partial f_1}{\partial x_2}(\bar x, \bar y) & \dots &  \frac{\partial f_1}{\partial x_n}(\bar x, \bar y)\\
		\frac{\partial f_2}{\partial x_1}(\bar x, \bar y) & \frac{\partial f_2}{\partial x_2}(\bar x, \bar y) & \dots &  \frac{\partial f_2}{\partial x_n}(\bar x, \bar y)\\
		\vdots & \vdots & \vdots & \vdots\\
		\frac{\partial f_p}{\partial x_1}(\bar x, \bar y) & \frac{\partial f_p}{\partial x_2}(\bar x, \bar y) & \dots &  \frac{\partial f_p}{\partial x_n}(\bar x, \bar y)\\
	\end{array}
	\right).
\]
Thus, for a scalar function $f:\mbbR^n\times\mbbR^m\rightarrow\mbbR$, $\nabla_xf(\bar x,\bar y)\in\mbbR^{1 \times p}$ is a row vector.
Given a function $f:\mbbR^n\times\mbbR^m\rightarrow\mbbR^p$ and an index set $\mathbb{I} = \{i_1,i_2,\dots,i_q\} \subseteq \{1,2,\dots,p\}$, the function $f_{\mathbb{I}}:\mbbR^n\times\mbbR^m\rightarrow\mbbR^q$, $q := \mathrm{dim}(\mathbb{I})$, is given by $f_{\mathbb{I}} = (f_{i_1}, f_{i_2}, \dots, f_{i_q})$. That is, $f_{\mathbb{I}}$ only has the $q$ elements of $f$ specified in $\mathbb{I}$.
The $n$-dimensional identity matrix is denoted by $\mathbf{I}_n\in\mbbR^{n\times n}$, and the matrix with $n$ rows and $m$ columns with only zero entries is denoted by $\mathbf{0}_{n \times m}\in\mbbR^{n\times m}$.
A vector of zero entries is denoted by $\zero$ (its dimension can always be inferred from context). 
Given a vector $x\in\mbbR^n$, $\Vert x \Vert$ denotes the Euclidean norm.  
For $x\in\mbbR^n$, $x \geq 0$ is interpreted element-wise, that is $x_i \geq 0$ for all $i = 1,2,\dots, n$.

\section{Preliminaries: Globally convergent homotopies for Systems of equations}\label{sec:homo_for_sys_eqns}

\begin{figure}[htb]
	\centering
	\includegraphics[scale=0.18]{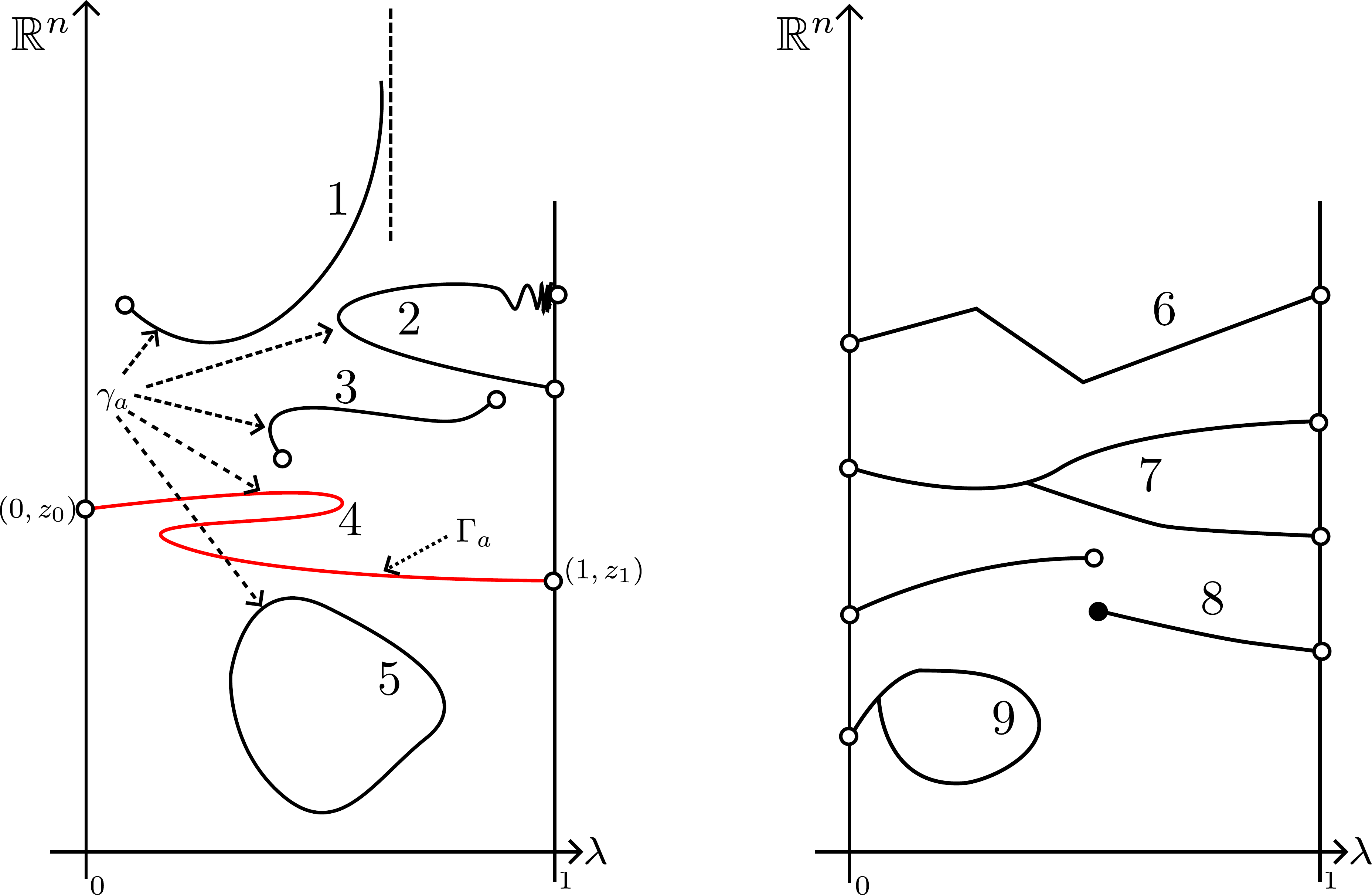}
	\caption{Under the conditions of Theorem~\ref{thm:trans} components of $\gamma_a$ are $C^1$ curves that are diffeomorphic to an interval (curves 1-4) or a circle (curve 5) and that may be unbounded (curve 1); have infinite arc length (curve 2); have (resp. not have) limit points $(0,z_0)$ and $(1,z_1)$ (curve 3, resp. curve 4). 
	They \emph{cannot} be nondifferentiable (curve 6); bifurcate (curve 7); be discontinuous (curve 8); or self-intersect (curve 9).}
	\label{fig:zero_curves}
\end{figure}

This section summarises some relevant theory on globally convergent homotopies for solving systems of equations, as introduced in~\cite{chow1978finding}. 
The reader may also consult the references~\cite{watson1979solving,watson2001theory}. 
Because this forms the fundamental basis on which the rest of the paper is built, we present this theory in some detail so that the paper is self-contained.
\begin{definition}
    Let $U\subset \mbbR^m$ and $V\subset \mbbR^n$ be open sets with $ m > n$. Then, $\rho \in \mathcal{C}^2(U,V)$ is said to be \emph{transversal to zero} if and only if
    \[
	\rank[\nabla \rho(z)] = n \qquad \forall\,z\in \rho^{-1}(\zero),
    \]
    where $\rho^{-1}(\zero) := \{z\in U : \rho(z) = \zero \}$.
\end{definition}
Consider now a $C^2$-map $\rho:\mbbR^m \times (0,1) \times \mbbR^n \rightarrow \mbbR^n$ and define $\rho_a:(0,1) \times \mbbR^n \rightarrow \mbbR^n$ as 
\[
    \rho_a(\lambda,z):= \rho(a,\lambda,z) \qquad\text{ with }\qquad  a \in \mbbR^m.
\]
Theorem~\ref{thm:trans} is a parametrised Sard's theorem, which is the main theoretical tool used throughout globally-convergent, probability-one homotopy theory. 
The proof can be found in~\cite{abraham1967transversal} and~\cite{chow1978finding}. 
\begin{theorem}\label{thm:trans}
    Let $\rho:\mbbR^m \times (0,1) \times \mbbR^n \rightarrow \mbbR^n$ be a $C^2$-mapping, which is transversal to zero. Then the mapping $\rho_a:(0,1) \times \mbbR^n \rightarrow \mbbR^n$ is transversal to zero for almost every $a\in\mbbR^m$ with respect to the Lebesgue measure.
\end{theorem}
As will be made clear in the sequel, the parameter $a$ will correspond to an initial point from where tracking of the zero curve, 
\[
\gamma_a = \{(\lambda, z) \in(0,1)\times \mbbR^n: \rho_a(\lambda, z) = \zero\}
\]
commences. 
The mapping $\rho_a$ being transversal to zero for almost every $a\in\mbbR^m$ means, by the implicit function theorem, that for almost every $a\in\mbbR^m$ the set $\gamma_a$ is a one-dimensional $C^1$-manifold. 
In other words, Theorem~\ref{thm:trans} says that if we have a $\rho$ 
transversal to zero and 
choose a random $a\in\mbbR^m$, then with probability one, $\gamma_a$ consists of a number of curves in $(0,1)\times\mbbR^{n}$ that are each diffeomorphic to a circle or an open interval, see also \cite[Cor.~2.3]{chow1978finding} and the discussions following this Corollary.  
Figure~\ref{fig:zero_curves} shows possible and impossible components of $\gamma_a$ (left and right curves, respectively).  

Next, we are interested in 
a component curve of~$\gamma_a$ (like the red curve in Figure~\ref{fig:zero_curves}), which is diffeomorphic to an interval, with points $(0,z_0)$ and $(1,z_1)$ as limit points for some well-designed map~$\rho$ 
so that $(0,z_0)$ is easily-obtained and the point $(1,z_1)$ corresponds to a solution of the difficult original problem.
The existence of such a component curve is guaranteed by Lemma~\ref{lem:homo_is_nice}, which is a slight adaptation of \cite[Lemmata~2.2 and~2.3]{watson2001theory}, see also~\cite{watson1979solving} as well as the proof of \cite[Thm.~2.5]{chow1978finding} for a version 
tailored to fixed-point problems. 
\begin{lemma}\label{lem:homo_is_nice}
	Consider a $C^2$-mapping $\rho:\mbbR^m\times[0,1)\times\mbbR^n$, which is transversal to zero when restricted to the open domain $\mbbR^m\times(0,1)\times\mbbR^n$. 
    Suppose that, for each $a\in\mbbR^m$, the system $\rho_a(0,z) = \zero$ has a unique nonsingular solution that is, there exists a $z_0\in\mbbR^n$ such that,
	\[
	\rho_a(0,z_0) = \zero \quad \text{ and } \quad \mathrm{rank}[\nabla \rho_a(0, z_0)] = n.
	\]
	Then, for almost all $a\in\mbbR^m$ there exists a unique $C^1$-curve, $\Gamma_a\subset \gamma_a$, emanating from $(0, z_0)$. 
	Moreover, if $\Gamma_a$ is bounded then it has a limit point at $(1,z_1)$, $z_1 \in \mbbR^m$ and if $\mathrm{rank}[\lim_{\lambda\rightarrow 1^-}\nabla \rho_a(\lambda, z_1)] = n$, then $\Gamma_a$ has finite arc length. 
\end{lemma}
The idea now is to construct maps (homotopies) that satisfy the conditions of Theorem~\ref{thm:trans} and Lemma~\ref{lem:homo_is_nice}. 
Then, for almost every $a\in\mbbR^m$ there exists a convenient curve, $\Gamma_a$, which one can track with numerical methods starting from the easily-obtained point $(0,z_0)$, until one arrives (possibly only arbitrarily closely) at the desired solution, $(1,z_1)$. 
Hence the name: globally-convergent probability-one homotopies.

\section{Globally convergent homotopies: Nonlinear programs}\label{sec:homoh_for_NLP}

In this section we summarise content from the paper \cite{watson2001theory} and show how the theory from Section~\ref{sec:homo_for_sys_eqns} can be used to solve nonlinear programs. 
We first explain how one can derive a suitable homotopy from the Karush-Kuhn-Tucker (KKT) conditions (denoted $\alpha$) and then present the paper's first contribution in Subsection~\ref{subsec:new_conditions}: a new assumption that renders the modified homotopy (denoted $\rho$) globally convergent. 

\subsection{From the KKT conditions to a system of equations}

\sloppy
Consider a parametrised nonlinear program ($\NLP^{\lambda}$), where $\lambda\in[0,1]$ is the homotopy parameter,
\begin{equation}\label{eq:NLP}
\NLP^{\lambda}:\begin{cases}
	\min\limits_{\boldu\in\mathbb{R}^r} & \quad J(\boldu),\\
	\mathrm{subject\  to:} & \quad G(\lambda, \boldu) \leq \zero.
\end{cases}
\end{equation}
Here $\boldu\in\mbbR^r$ is the decision variable, $J:\mbbR^r\rightarrow \mbbR$ is the cost function and $G~:~[0,1]~\times~\mbbR^r\rightarrow \mbbR^s$ denotes the parametrised constraints. 
Throughout the paper we assume:
\begin{enumerate}
	\item[\textbf{(A1)}] Both $J$ and $G$ are $C^3$ functions. 
\end{enumerate}
We assume that with $\lambda = 0$ the problem is easy to solve, whereas the original difficult problem is obtained when $\lambda = 1$. 
The motivation for perturbing the constraints comes from our intended application of path planning problems for dynamical systems where obstacles often appear as constraints that are nonconvex with respect to the state. As we show in our numerics section, we will introduce $\lambda$ such the obstacles ``shrink'' to points when $\lambda = 0$, and ``grow'' to their original sizes as $\lambda$ increases to 1. 

\fussy
The KKT conditions, see for example \cite[Ch. 5]{boyd2004convex} \cite[Ch. 3]{bertsekas1999nonlinear}, state that if a constraint qualification holds and if $\boldu^{\lambda}\in\mbbR^r$ locally solves $\NLP^{\lambda}$, then there exists a $\boldsymbol \mu^{\lambda}\in\mbbR^s$ such that:
\begin{align}
	(\nabla_\boldu J(\boldu^{\lambda}))^\top +  \left(\nabla_\boldu G(\lambda, \boldu^{\lambda})\right)^\top \boldsymbol \mu^{\lambda} &= \zero,\label{eq:KKT0}\\
	G(\lambda,\boldu^{\lambda}) & \leq \zero,\label{eq:KKT1}\\
	(\boldsymbol \mu^{\lambda})^\top G(\lambda, \boldu^{\lambda}) &= 0,\label{eq:KKT2}\\
	\boldsymbol \mu^{\lambda} & \geq \zero.\label{eq:KKT3}
\end{align}
Throughout the rest of the paper we will let $(\boldu^\star, \boldsymbol \mu^\star)$ denote a primal-dual pair satisfying \eqref{eq:KKT0}-\eqref{eq:KKT3} with $\lambda = 1$, that is, $\boldu^\star = \boldu^1$ and $\boldsymbol \mu^\star = \boldsymbol \mu^1$.

There are many constraint qualifications, some of which are stricter than others, see for example \cite[Ch. 3]{bertsekas1999nonlinear}. 
For example, a common one is the \emph{linear independence constraint qualification} (LICQ), which states that the gradients of the constraints that are active at $\boldu^\lambda$ must be linearly independent. That is, the set $\{\nabla_{\boldu}G_i(\lambda, \boldu^{\lambda}) : i\in\mathbb{I}(\boldu^{\lambda}) \}$ must be linearly independent, $\mathbb{I}(\boldu^{\lambda}) := \{i\in\{1,\dots,s\}: G_i(\lambda, \boldu^{\lambda}) = \zero\}$. 

The conditions \eqref{eq:KKT1}-\eqref{eq:KKT3} by themselves form a \emph{complementarity problem}, which can be written as a system of equations via the following theorem due to Mangasarian, see \cite{mangasarian1976equivalence}. 
\begin{theorem}\label{thm:Mang}
	Consider the complementarity problem: find $z\in\mbbR^n$ such that $z^\top F(z) = 0$, $F(z) \geq \zero$, $z\geq \zero$, where $F:\mbbR^n\rightarrow \mbbR^n$. 	
	Let $\theta:\mbbR\rightarrow\mbbR$ be any strictly increasing function and let $\theta(0) = 0$. 
	Then $z$ solves the complementarity problem if and only if
	\[
		\theta(| F_i(z) - z_i|) - \theta(F_i(z)) - \theta(z_i) = 0,\quad \forall\, i\in[1:n].
	\]
\end{theorem}
Following \cite{watson1979solving}, the function $\theta(t) = t^3$ is a simple choice for which $\theta(|t|)$ is $C^2$. 
(Referring back to Theorem~\ref{thm:trans}, we will need a $C^2$ mapping in the sequel).  
Thus, if we define,
\begin{equation}
	\alpha(\lambda,\boldu,\boldsymbol \mu):=
	\left[
	\begin{array}{c}
		( \nabla_{\boldu} J(\boldu) )^\top + \left(\nabla_{\boldu} G(\lambda, \boldu)\right)^\top \boldsymbol \mu \\
		(\mu_1)^3 - | - G_1(\lambda, \boldu) - \mu_1|^3 + (-G_1(\lambda, \boldu))^3 \\
		(\mu_2)^3 - | - G_2(\lambda, \boldu) - \mu_2|^3 + (-G_2(\lambda, \boldu))^3\\
		\vdots\\
		(\mu_s)^3 - | - G_s(\lambda, \boldu) - \mu_s|^3 + (-G_s(\lambda, \boldu))^3 
	\end{array}
	\right],\label{eq:alpha}
\end{equation}
where $\boldsymbol \mu := (\mu_1,\mu_2,\dots, \mu_s)^\top$, then, via Theorem~\ref{thm:Mang}, for a given $\lambda\in[0,1]$, the pair $(\boldu^{\lambda}, \boldsymbol \mu^{\lambda})$ satisfies the KKT conditions \eqref{eq:KKT0}-\eqref{eq:KKT3} if and only if
\[
\alpha(\lambda,\boldu^\lambda,\boldsymbol \mu^\lambda) = \boldsymbol 0.
\]

\subsection{Modifications to $\alpha$ to obtain a globally convergent homotopy}

We would like to traverse the set $\alpha^{-1}(\boldsymbol 0) = \{(\lambda,\boldu,\boldsymbol \mu)\in[0,1]\times\mbbR^r\times\mbbR^s : \alpha(\lambda,\boldu,\boldsymbol \mu) = \boldsymbol 0\}$, starting from the easily obtained point, $(0,\boldu^0,\boldsymbol \mu^0)$, towards a point $(1,\boldu^\star,\boldsymbol \mu^\star)$ where we arrive at a local solution to $\NLP^{1}$. 
However, the components of $\alpha^{-1}$ may not be convenient to track numerically. 
We now describe how $\alpha$ can be modified to obtain a homotopy, $\rho$, which meets the conditions of Theorem~\ref{thm:trans} and Lemma~\ref{lem:homo_is_nice}, so that it is globally convergent with probability one. 

Consider an arbitrary nonempty open set,
\begin{equation}
	U^0\subset\mbbR^r, \label{eq:U0}
\end{equation}
and define,
\begin{equation}
	B^0(\boldu) := \{b\in\mbbR^s_{> 0}: G(0, \boldu) < b  \},\quad \boldu \in U^0. \label{eq:B0}
\end{equation}
Thus, for some $\boldu\in U^0$, the set $B^0(\boldu)$ contains all positive vectors that say how much the constraints can be relaxed such that $\boldu$ is feasible when $\lambda = 0$. 
Now, define the homotopy, $\rho:\mbbR^{r} \times \mbbR^s_{>0} \times \mbbR^s_{>0}\times [0,1)\times\mbbR^r\times\mbbR^s\rightarrow \mbbR^{r + s}$,
\begin{equation}
	\rho(\boldu^0, b^0, c^0, \lambda, \boldu, \boldsymbol \mu) := 
	\left(
	\begin{array}{c}
		\lambda [ (\nabla_{\boldu} J(\boldu))^\top + \left(\nabla_{\boldu} G(\lambda,\boldu)\right)^\top \boldsymbol \mu] + (1 - \lambda)  (\boldu- \boldu^0) \\
		K(\lambda, \boldu, \boldsymbol \mu, b^0, c^0)
	\end{array}
	\right), \label{eq_homo}
\end{equation}
where $K:[0,1)\times \mbbR^r \times \mbbR^s \times \mbbR^s_{>0} \times \mbbR^s_{>0} \rightarrow \mbbR^s$ is defined as follows, for $i\in[1:s]$,
\begin{equation}
	K_i(\lambda, \boldu, \boldsymbol \mu, b^0, c^0) := \mu_i^3 - |(1 - \lambda)b_i^0 - G_i(\lambda, \boldu) - \mu_i |^3 + \left[ (1 - \lambda)b_i^0   - G_i(\lambda, \boldu)\right]^3  - (1 - \lambda)c_i^0, \label{eq:K}
\end{equation}
where $b^0 := (b_1^0, b_2^0,\dots, b_s^0)^\top$ and $c^0 := (c_1^0, c_2^0,\dots, c_s^0)^\top$.
Finally, we introduce the \emph{parametrised homotopy},
\begin{equation}
\rho_a(\lambda, \boldu, \boldsymbol \mu) = \rho(\boldu^0, b^0, c^0, \lambda, \boldu, \boldsymbol \mu), \label{def:para_homo}
\end{equation}
where $a = (\boldu^0,b^0,c^0)\in U^0 \times B^0(\boldu^0)\times \mbbR^s_{>0}$ is an arbitrarily chosen point, and let,
\[
\gamma_a = \rho_a^{-1}( \boldsymbol 0) := \{(\lambda, \boldu, \boldsymbol \mu)\in[0,1)\times\mbbR^r\times\mbbR^s: \rho_a(\lambda, \boldu, \boldsymbol \mu) = \boldsymbol 0\}. 
\]
When compared with $\alpha$, see \eqref{eq:alpha}, three new parameters appear in $\rho$.
The point $\boldu^0\in U^0$ is an \emph{arbitrary} initial decision vector. 
Note that this point does not even need to be a feasible point of the relaxed problem. 
As we will show in the proof of Proposition~\ref{prop:main:NLPs}, the vector $c^0\in\mbbR^s_{>0}$ is introduced so that the rank of the matrix $\nabla_{c^0} \rho$ is full, which is needed to render $\rho$, see \eqref{eq_homo}, transversal to zero so that we may invoke Theorem~\ref{thm:trans}.  
The vector $b^0\in\mbbR^s_{>0}$ is introduced to keep the (what will be shown to be unique) curve $\Gamma_a\subset\gamma_a$ bounded so that we may invoke Lemma~\ref{lem:homo_is_nice} to deduce that $\Gamma_a$ reaches a point $(1,\boldu^\star, \boldsymbol \mu^\star)$.  

An important point to appreciate is that, unlike with the homotopy $\alpha$, a point $(\lambda, \boldu^{\lambda}, \boldsymbol \mu^{\lambda})$ for which $\rho_a(\lambda, \boldu^{\lambda}, \boldsymbol \mu^{\lambda}) = \zero$, with $\lambda \neq 1$, does \emph{not} correspond to a solution of the KKT conditions, \eqref{eq:KKT0}-\eqref{eq:KKT3}.
Recall that the idea is to easily find the point $(0,\boldu^0,\boldsymbol \mu^0)$, and to know that there exists a ``nice'' path to track towards $(1,\boldu^\star, \boldsymbol \mu^\star)$, the desired solution. 

We conclude this section with a lemma that collects a few properties of the function $K$, stated throughout the paper \cite{watson2001theory}. 
This will make the proofs of Proposition~\ref{prop:main:NLPs} and Proposition~\ref{prop:main_result_OCP} easier to digest.
\begin{lemma}\label{lem:K_properties}
	For an arbitrary $i\in[1:s]$, consider the function $K_i$, as defined in \eqref{eq:K}, with $\lambda\in[0,1)$, $b^0\in\mbbR^s_{>0}$ and $c^0\in\mbbR^s_{>0}$.
	The following holds.
	\begin{enumerate}
		\item[1.]  If $\boldu\in\mbbR^r$ is chosen such that $G_i(\lambda,\boldu) < (1-\lambda)b_i^0$, then $K_i$ is a strictly increasing function of $\mu_i$.
		\item[2.] If $\mu_i < 0$ or $G_i(\lambda, \boldu) > (1 - \lambda)b_i^0$, then $K_i < 0$.
		\item[3.] For a given $a = (\boldu^0, b^0, c^0)\in U^0\times B^0(\boldu^0)\times \mbbR^s_{>0}$, there exists a unique $\boldsymbol \mu^0\in\mbbR^s_{> 0}$ such that $K(0,\boldu^0,\boldsymbol \mu^0,b^0,c^0) = \boldsymbol 0$.
	\end{enumerate}
\end{lemma}
\begin{proof}
	\textbf{Point~1.} If $(1 - \lambda)b_i^0 - G_i(\lambda, \boldu) - \mu_i  <  0$. 
	Then,
	\[
	\frac{\partial K_i}{\partial \mu_i} =  3\mu_i^2 - 3( - (1 - \lambda)b_i^0 + G_i(\lambda, \boldu) + \mu_i)^2.
	\]
	Recall we assume that $G_i(\lambda,\boldu) < (1-\lambda)b_i^0$, therefore,
	\[
 \mu_i^2 \geq \left(-(1 - \lambda)b_i^0 + G_i(\lambda,\boldu) + \mu_i\right)^2.
	\]
	Therefore,
	\[
	3\left( \mu_i^2 - \left(-(1 - \lambda)b_i^0 + G_i(\lambda,\boldu) + \mu_i\right)^2  \right) = \frac{\partial K_i}{\partial \mu_i}\geq  0.
	\]
	Moreover,
	$\frac{\partial K_i}{\partial \mu_i}\neq  0$, because then $G_i(\lambda,\boldu) = (1-\lambda)b_i^0$, contradicting our assumption that $G_i(\lambda,\boldu) < (1-\lambda)b_i^0$.
	Thus, $\frac{\partial K_i}{\partial \mu_i}>  0$.
	
	If $(1 - \lambda)b_i^0 - G_i(\lambda, \boldu) - \boldsymbol \mu_i \geq 0$ we immediately get,
	\[
	\frac{\partial K_i}{\partial \mu_i} =  \mu_i^2 + 3((1 - \lambda)b_i^0 - G_i(\lambda, \boldu) - \mu_i)^2 \geq 0,
	\]
	and the argument that $\frac{\partial K_i}{\partial \mu_i} \neq  0$ is the same as before.
	
	\noindent\textbf{Point~2.} If  $\mu_i < 0$, then 
	\[
	(1-\lambda)b_i^0 - G_i(\lambda,\boldu) - \mu_i \geq  (1-\lambda)b_i^0 - G_i(\lambda,\boldu).
	\] 
	Therefore,
	\[
	|(1-\lambda)b_i^0 - G_i(\lambda,\boldu) - \mu_i|^3  \geq [(1-\lambda)b_i^0 - G_i(\lambda,\boldu)]^3.
	\]
	We then see that:
	\[
	K_i = \underbrace{ \mu_i^3}_{< 0} \underbrace{- |(1-\lambda)b_i^0 - G_i(\lambda,\boldu) -  \mu_i|^3  + [(1-\lambda)b_i^0 - G_i(\lambda,\boldu)]^3 }_{\leq 0} \underbrace{- (1-\lambda)c_i^0}_{< 0} < 0.
	\]
	If $G_i(\lambda, \boldu) > (1 - \lambda)b_i^0$ and $\mu_i < 0$, then we immediately get:
	\[
	K_i = \underbrace{ \mu_i^3}_{< 0} \underbrace{- |(1-\lambda)b_i^0 - G_i(\lambda,\boldu) -  \mu_i|^3}_{<0}  + \underbrace{[(1-\lambda)b_i^0 - G_i(\lambda,\boldu)]^3 }_{< 0} \underbrace{- (1-\lambda)c_i^0}_{< 0} < 0,
	\]
	If $G_i(\lambda, \boldu) > (1 - \lambda)b_i^0$ and $ \mu_i \geq 0$, then
	\[
	(1 - \lambda)b_i^0 - G_i(\lambda, \boldu) -  \mu_i < - \mu_i \leq 0.
	\]
	Moreover,
	\[
	|(1 - \lambda)b_i^0 - G_i(\lambda, \boldu) -  \mu_i | > | -  \mu_i | =  \mu_i.
	\]
	Therefore,
	\[
	 \mu_i^3 < |(1 - \lambda)b_i^0 - G_i(\lambda, \boldu) -  \mu_i |^3.
	\]
	Thus,
	\[
	K_i = \underbrace{ \mu_i^3 - |(1-\lambda)b_i^0 - G_i(\lambda,\boldu) -  \mu_i|^3}_{<0}  + \underbrace{[(1-\lambda)b_i^0 - G_i(\lambda,\boldu)]^3 }_{< 0} \underbrace{- (1-\lambda)c_i^0}_{< 0} < 0.
	\]
	\textbf{Point~3.} From the definition of $B^0(\boldu)$, see~\eqref{eq:B0}, we get $G_i(0,\boldu^0) - b_i^0 < 0$.
	Thus,
	\[
	K_i(0,\boldu^0,\mathbf{0}_s,b^0,c^0) = - |b_i^0 - G_i(0, \boldu^0)|^3 + \left[ b_i^0   - G_i(0, \boldu^0)\right]^3  - c_i^0 = -c_i^0 < 0.
	\]
	Moreover, for every $i$ there exists a finite $\bar{\boldsymbol \mu}_i\in\mbbR_{> 0}$ large enough such that,
	\[
	\underbrace{\bar{\mu}_i^3}_{>0} \underbrace{- |b_i^0 - G_i(0,\boldu^0) - \bar{ \mu}_i|^3 -c_i^0 }_{>0}  + \underbrace{[b_i^0 - G_i(0,\boldu^0)]^3 }_{> 0} > 0.
	\]
	Thus, from the intermediate value theorem and the fact that $K_i$ is strictly increasing with respect to $ \mu_i$ (see Point 1.) we can conclude that there exists a unique $ \mu_i^0$ satisfying $0< \mu_i^0<\bar{ \mu}_i$ for which
	\[
	( \mu_i^0)^3- |b_i^0 - G_i(0,\boldu^0) -  \mu_i^0|^3 -c_i^0  + [b_i^0 - G_i(0,\boldu^0)]^3 = 0.
	\]
	The result then follows with $\boldsymbol \mu^0\in\mbbR^s_{>0}$ defined as $ \boldsymbol\mu^0 := ( \mu_1^0, \mu_2^0,\dots,  \mu_s^0)^\top$. 
\end{proof}

\subsection{A new condition to render $\rho$ globally convergent}\label{subsec:new_conditions}

Up to this point we have mostly summarised previously know results on convergent homotopies. 
We now present the paper's first contribution, which says that $\rho$ can be rendered globally convergent if, along with a few other assumptions, the homotopy parameter is introduced such that the nonconvex constraints are trivially satisfied with $\lambda = 0$ (Assumption~(A2)). 

Referring to $\NLP^{\lambda}$, see \eqref{eq:NLP}, we distinguish between constraints that are convex and nonconvex with respect to $\boldu$ by introducing,
\[
\mathbb{I}_c := \{i\in\{1,2,\dots, s\}: G_i(\lambda,\boldu)\ \mathrm{is\ convex\ w.r.t\ }\boldu,\  \forall \lambda\in[0,1]\ \},
\]
and
\[
\mathbb{I}_{nc} := \{i\in\{1,2,\dots,s\}\setminus\mathbb{I}_c\}.
\]
The new assumption we impose reads as follows,
\begin{enumerate}
	\item[\textbf{(A2)}] $G_i(0, \boldu) \leq 0$ for all $\boldu\in\mbbR^r$, for all $i \in \mathbb{I}_{\mathrm{nc}}$.
\end{enumerate} 
This states that all nonconvex constraints are trivially satisfied for the totally relaxed problem, where $\lambda = 0$. 
We emphasise that (A2) needs to hold for \emph{all} $\boldu\in\mbbR^r$, not just for an initial $\boldu^0$ corresponding to $a$, and not just for all $\boldu$ for which the convex constraints are satisfied. 
As we will show, (A2) implies that once we start traversing the zero curve, $\Gamma_a$, we cannot loop back to $\lambda = 0$. 

The next assumption, which appears in \cite[Thm 7.1]{watson2001theory}, concerns the set,
\begin{equation}
	S(\lambda, b) := \{\boldu\in\mbbR^r : G(\lambda, \boldu) \leq (1 - \lambda) b \}, \quad \lambda\in[0,1],\,\, b \in\mbbR^s_{> 0},\label{def:S}
\end{equation}
and reads as follows,
\begin{enumerate}
	\item[\textbf{(A3)}] For a chosen pair, $\boldu^0\in U^0$  and $b^0\in B^0(\boldu^0)$, we have
	\[
		S(\lambda, b^0) \neq \emptyset,\quad  \forall \lambda \in[0,1] \quad \text{and} \quad \cup_{\lambda\in[0,1]} S(\lambda , b^0)\quad \text{is bounded}.
	\]
\end{enumerate}
If $S(\lambda,b^0)$ is empty for some $\lambda\in[0,1]$, then there exists an $i$ such that $G_i(\lambda,\boldu) > (1 - \lambda)b_i^0$, thus $K_i < 0$ (see Point.~2 from Lemma~\ref{lem:K_properties}) and thus $\rho_a(\lambda,\boldu,\mu) \neq \zero$. 
Thus, we require $S(\lambda,b^0)$ nonempty for all $\lambda\in[0,1]$. 
Boundedness of $\cup_{\lambda\in[0,1]} S(\lambda , b^0)$ will be needed to argue that $\Gamma_a$ is bounded, allowing us to invoke Lemma~\ref{lem:homo_is_nice}. 

The next assumption, also taken directly from \cite[Thm 7.1]{watson2001theory}, may be interpreted as a constraint qualification for homotopies. 
\begin{itemize}
	\item[\textbf{(A4)}] For any limit point $(\hat \lambda, \hat \boldu)$ along $\gamma_a$ there exists a $d\in\mbbR^r$ such that
	\[
	\nabla_{\boldu} G_{\hat{\mathbb{I}}}(\hat \lambda, \hat \boldu)d > \zero,\  \mathrm{where}\  \hat{\mathbb{I}} :=  \{i\in\{1,2,\dots,s\} :G_i(\hat \lambda, \hat \boldu) =  (1 - \hat \lambda)b_i^0\}.
	\]
\end{itemize}
We of course need the following assumption to hold at the point $(1,\boldu^{\star}, \boldsymbol \mu^{\star})$ so that the KKT conditions holds.
\begin{itemize}
	\item[\textbf{(A5)}] A constraint qualification holds at every local solution of $\NLP^1$. 
\end{itemize}
Finally, we explicitly write out the following assumption, which appears in Theorem~\ref{thm:trans}, ensuring that $\Gamma_a\subset\gamma_a$ has finite arc length.
\begin{itemize}
	\item[\textbf{(A6)}] The $\mathsf{rank}\left[ \nabla \rho_a(1, \boldu^{\star} , \boldsymbol  \mu^{\star}) \right] = r + s$.
\end{itemize}

As we will show in the next section, Assumptions (A2)-(A3) are concerned with how the homotopy parameter is introduced. 
The Assumptions (A4)-(A6) may be considered as ``technical assumptions'', which are difficult to verify beforehand, but only fail in special cases. 

We now state our first result.
\sloppy
\begin{proposition}\label{prop:main:NLPs}
	Consider $\NLP^\lambda$, with $\lambda\in[0,1]$, as in \eqref{eq:NLP}, under Assumption~(A1), along with an arbitrary open set $U^0\subset \mbbR^r$, see \eqref{eq:U0}, the set $B^0(\boldu)$, see \eqref{eq:B0}, and the homotopy $\rho$, see \eqref{eq_homo}. 
	The following holds: for almost every $a = (\boldu^0, b^0, c^0) \in U^0 \times B^0(\boldu^0) \times \mbbR^s_{> 0}$ there exists a unique curve, $\Gamma_a\subset \gamma_a$, emanating from the unique point $(0, \boldu^0, \boldsymbol \mu^0)$, which is $C^1$, does not intersect itself or other components of $\gamma_a$, and does not bifurcate. (See Figure~\ref{fig:zero_curves} in Section~\ref{sec:homo_for_sys_eqns}). 	
	
	Furthermore, suppose that Assumptions~(A2)-(A5) hold as well. 
	Then $\Gamma_a$ has a limit point, $(1,\boldu^\star, \boldsymbol \mu^\star)$, where the pair $(\boldu^\star, \boldsymbol \mu^\star)$ is a local solution to $\NLP^1$. 
	If, additionally, (A6) holds, then $\Gamma_a$ has finite arc length. 
\end{proposition}
\fussy
\begin{proof}
	The proof collects a number of arguments made in the various proofs and discussions throughout the paper \cite{watson2001theory}. 
	Because $J$ and $G$ are $C^3$, $\rho$ is $C^2$. 
	The Jacobian matrix $\nabla \rho$ has rank $r + s$ for all $\lambda \in [0,1)$, because the matrices $\nabla_{\boldu^0}\rho$ and $\nabla_{c^0}\rho$ read:
	\[
	\nabla_{\boldu^0}\rho = -(1-\lambda)\left(
	\begin{array}{c}
		\mathbf{I}_r\\
		\mathbf{0}_{s\times r}
	\end{array}
	\right)
	\quad \text{and} \quad
		\nabla_{c^0}\rho = -(1-\lambda)
	\left(
	\begin{array}{c}
		\mathbf{0}_{r \times s}\\
		\mathbf{I}_{s}
	\end{array}
	\right),
	\]
	which have rank $r$ and $s$, respectively. 
	Thus, $\rho$ is transversal to zero on $\mbbR^{r} \times \mbbR^s_{>0} \times \mbbR^s_{>0}\times (0,1)\times\mbbR^r\times\mbbR^s$, and we may invoke Theorem~\ref{thm:trans} to immediately deduce that  $\rho_a$ is transversal to zero for almost every $a=(\boldu^0, b^0, c^0)\in U^0 \times B^0(\boldu^0) \times \mbbR^s_{>0}$. 
	By Point~3. of Lemma~\ref{lem:K_properties}, for any $a$ there exists a unique point $\boldsymbol \mu^0\in\mbbR^s_{> 0}$ such that $\rho_a(0,\boldu^0,\boldsymbol \mu^0) = \zero$. 
	We can thus invoke Lemma~\ref{lem:homo_is_nice} to deduce the existence of the unique curve $\Gamma_a$ with the properties listed in this proposition. 
	
	We now argue that (A2)-(A5) imply that $\Gamma_a$ is bounded to conclude, via the second part of Lemma~\ref{lem:homo_is_nice}, that $\Gamma_a$ reaches the point $(1,\boldu^{\star}, \boldsymbol \mu^{\star})$.
	This part of the proof follows the same lines of reasoning as the proof of \cite[Thm 7.1]{watson2001theory}, the only difference being that we substitute a key assumption (which prevents $\Gamma_a$ from looping back to $\lambda = 0$ as one traverses along it) with the new assumption, (A2).
	We will present those details of the proof affected by the new assumption and only sketch those parts of the proof of \cite[Thm 7.1]{watson2001theory} not affected by it.
	
	Consider an arbitrary point $(\lambda, \boldu,\boldsymbol \mu) \in \Gamma_a$, with $\lambda\in[0,1)$. 
	Recall that for such a point $K_i = 0$ for all $i\in[1:s]$ and thus, by~Point 2. of Lemma~\ref{lem:K_properties}, we have $\mu_i \geq 0$ and $G_i(\lambda, \boldu) \leq (1-\lambda) b_i^0$ for all $i\in[1:s]$.
	Thus, if $(\lambda, \boldu,\boldsymbol \mu) \in \Gamma_a$, then $\boldu \in S(\lambda,b^0)$, $b^0\in\mbbR^s_{>0}$. 
	Seeking a contradiction, suppose now that $\Gamma_a$ is unbounded and let $\{(\lambda^k, \boldu^k, \boldsymbol \mu^k)\}_{k\in\mathbb{N}}\subset\Gamma_a$ be a sequence for which $\Vert (\lambda^k, \boldu^k, \boldsymbol \mu^k)\Vert \rightarrow \infty$. 
	Because $[0,1] \times \cup_{\lambda\in[0,1]} S(\lambda,b^0)$ is nonempty and compact (by (A3)) there exists a subsequence $\{(\lambda^{k_j}, \boldu^{k_j}, \boldsymbol \mu^{k_j})\}$, with $\boldsymbol \mu^{k_j} \in\mbbR^s_{\geq 0}$ such that $\{(\lambda^{k_j}, \boldu^{k_j})\} \rightarrow (\hat \lambda, \hat \boldu)\in [0,1] \times \cup_{\lambda\in[0,1]} S(\lambda,b^0)$, but with $\Vert \boldsymbol \mu^{k_j} \Vert \rightarrow \infty$. 
	
	From the homotopy, \eqref{eq_homo}, we have,
	\begin{align}
		\hat \lambda [(\nabla_\boldu J(\hat \boldu))^\top + \left(\nabla_\boldu G(\hat \lambda, \hat \boldu)\right)^\top \boldsymbol \mu^{k_j}] + (1 - \hat \lambda)(\hat \boldu- \boldu^0) & \rightarrow \zero,\label{eq_grad_limit}\\
		K_i(\hat \lambda, \hat \boldu, \boldsymbol \mu^{k_j}, b^0, c^0) & \rightarrow 0, \quad i=1,2,\dots,s. \label{eq_K_limit}
	\end{align}
	Recall that if $(\lambda,\boldu,\boldsymbol \mu)\in\Gamma_a$ then $\mu_i \geq 0$ and $G_i(\lambda,\boldu)\leq (1-\lambda)b_i^0$ for all $i\in[1:s]$. 
	Seeking a contradiction, suppose that there exists an $\tilde i\in[1:s]$ for which $G_{\tilde i}(\hat \lambda,\hat\boldu) <  (1-\hat \lambda)b_{\tilde i}^0$ and $\limsup_{k_j \rightarrow\infty} \mu_{\tilde i}^{k_j} = \infty$. 
	Then, from Point~1. of Lemma~\ref{lem:K_properties}, $K_{\tilde i}(\hat \lambda, \hat \boldu, \boldsymbol \mu^{k_j}, b^0, c^0)$ would be strictly increasing. 
	Moreover, recall from the proof of Point~3. of Lemma~\ref{lem:K_properties} that $K_i > 0$ if $\mu_i$ is large enough. 
	Therefore, referring to \eqref{eq_K_limit}, we would have a contradiction. 
	We may conclude that for every $i$ for which $\limsup_{k_j \rightarrow\infty} \mu_i^{k_j} = \infty$, we have $G_i(\hat \lambda,\hat\boldu) =  (1-\hat \lambda)b_i^0$.
	Denote the set of these indices by $\tilde{\mathbb{I}}$. That is,
	\[
	\tilde{\mathbb{I}} := \{i\in\{1,2,\dots,s\} : \limsup_{k_j \rightarrow\infty} \mu_i^{k_j} = \infty \}.
	\]
	Suppose now that $\hat \lambda = 0$. Then, $G_i(0, \hat \boldu) =  b_i^0 > 0$ for all $i\in\tilde{\mathbb{I}}$. 
	This function $G_i$ cannot be a nonconvex constraint, that is, $i\notin\mathbb{I}_{nc}\cap\tilde{\mathbb{I}}$ because then Assumption (A2) would be violated.
	Moreover, $i$ cannot correspond to a convex constraint, that is $i\notin\mathbb{I}_{c}\cap\tilde{\mathbb{I}}$, the argument being the same as the one in \cite[Thm. 6.1]{watson2001theory}, which is as follows. 
	First, let 
	\[
	\tilde{\mathbb{I}}_c := \mathbb{I}_c \cap \tilde{\mathbb{I}}.
	\]	
	If $\hat \lambda = 0$ then from \eqref{eq_grad_limit}, $\lambda^{k_j} (\boldsymbol \mu^{k_j})^\top \nabla_{\boldu}G(0, \hat \boldu) \rightarrow (\boldu^0 - \hat \boldu)^\top$. 
	Thus, there exists a $w\in\mbbR_{\geq 0}^{s}$, for which $\lambda^{k_j} (\boldsymbol \mu^{k_j})^\top \rightarrow w$, satisfying $w_i = 0$ if $i\notin \tilde{\mathbb{I}}$, such that $w^\top\nabla_{\boldu}G(0, \hat \boldu) = (\boldu^0 - \hat \boldu)^\top$, or equivalently,
	\begin{equation}
		w_{\tilde{\mathbb{I}}_c}^\top\nabla_{\boldu}G_{\tilde{\mathbb{I}}_c}(0, \hat \boldu) = (\boldu^0 - \hat \boldu)^\top.\label{eq:w_grad_1}
	\end{equation} 
	Recall that if $\hat \lambda = 0$ then $G_{\tilde{\mathbb{I}}}(0,\hat \boldu) = b_{\tilde{\mathbb{I}}}^0$, and from the definition of $B^0(\boldu)$, see~\eqref{eq:B0}, $G_{\tilde{\mathbb{I}}}(0,\boldu^0) < b^0_{\tilde{\mathbb{I}}}$ for $\boldu^0\in U^0$.
	Therefore, also using convexity of $G_{\tilde{\mathbb{I}}_c}$, we have,
	\begin{equation}
		\nabla_{\boldu}G_{\tilde{\mathbb{I}}_c}(0, \hat \boldu)(\boldu^0 - \hat\boldu) \leq G_{\tilde{\mathbb{I}}_c}(0, \boldu^0) - G_{\tilde{\mathbb{I}_c}}(0, \hat \boldu) < b^0_{\tilde{\mathbb{I}}_c} - b^0_{\tilde{\mathbb{I}}_c} = \zero,\label{eq:w_grad_2}
	\end{equation}
	and therefore, multiplying the left-hand side of \eqref{eq:w_grad_2} by $w_{\tilde{\mathbb{I}}_c}^\top$ we get,	
	\[
	w_{\tilde{\mathbb{I}}_c}^\top \nabla_{\boldu}G_{\tilde{\mathbb{I}}_c}(0, \hat \boldu)(\boldu^0 - \hat\boldu) < 0.
	\]
	However, multiplying \eqref{eq:w_grad_1} on the right by $(\boldu^0 - \hat \boldu)$, we have,
	\[
	w_{\tilde{\mathbb{I}}_c}^\top \nabla_{\boldu}G_{\tilde{\mathbb{I}}_c}(0, \hat \boldu)(\boldu^0 - \hat\boldu) = (\boldu^0 - \hat\boldu)^\top(\boldu^0 - \hat \boldu) \geq 0.
	\]
	So we have a contradiction, and thus $\hat \lambda \neq 0$ for any $i\in \tilde{\mathbb{I}}$.
	
	The remainder of the proof is exactly the same as in \cite[Thm. 7.1]{watson2001theory}, which uses a technique from the proof of \cite[Thm. 5.1]{watson2001theory}. 
	We sketch the remainder to show where (A4) comes in. 
	Because $\hat \lambda \neq 0$, and due to the boundedness of $S(\lambda, b)$ for $\lambda\in[0,1]$ and $b\in\mbbR^s_{>0}$, one can argue that there exists a $w_2\in\mbbR^s_{\geq0}$ such that,
	\[
	\left(\nabla_\boldu J(\hat \lambda, \hat \boldu)\right)^\top + \left(\nabla_\boldu G(\hat \lambda, \hat \boldu)\right)^\top w_2 + \frac{(1 - \hat \lambda)(\hat \boldu- \boldu^0)}{\hat \lambda} = \zero,
	\]
	this equation coming from \eqref{eq_grad_limit}. If we let $\boldsymbol \mu^{k_j} = w_2 + v^{k_j}$, then $\Vert v^{k_j} \Vert \rightarrow\infty$, and it can be argued that the sequence $\{ v^{k_j}/\Vert v^{k_j} \Vert _{\infty} \}$ is bounded. Thus there exists a subsequence, converging to a point $v\in\mbbR^s_{\geq 0}$ for which $\Vert v \Vert_{\infty} = 1$. This then leads to the system $\left(\nabla_\boldu G(\hat \lambda, \hat \boldu)\right)^\top v = \zero$, or,	
	\[
	\left(\nabla_\boldu G_{\hat{\mathbb{I}}}(\hat \lambda, \hat \boldu)\right)^\top v = \zero,\quad v\in\mbbR^s_{\geq 0}\quad v \neq \zero,
	\]
	$\hat{\mathbb{I}}$ as defined in Assumption~(A4).
	Then, by Gordan's theorem, see for example \cite[App~A.4]{wright2022optimization}, there does \emph{not exist} a $d\in\mbbR^r$ such that,
	\[
	\nabla_\boldu G_{\hat{\mathbb{I}}}(\hat \lambda, \hat \boldu)d > 0, \quad d\in\mbbR^r
	\]
	which violates Assumption (A4). This is a contradiction, and so $\Gamma_a$ is bounded.
	By Lemma~\ref{lem:homo_is_nice}, we can conclude that $\Gamma_a$ has the stated limit point.
	The final statement involving (A6) also directly follows from Lemma~\ref{lem:homo_is_nice}.
\end{proof}

\section{Globally convergent homotopies: Discrete-time optimal control}\label{sec:homo_for_OCPs}

\sloppy
When wanting to solve an OCP via the globally convergent homotopy, \eqref{eq_homo}, it might not always be obvious how one should introduce the homotopy parameter, $\lambda$, such that (A2) and (A3) hold. 
This is due to the presence of dynamics and the fact that the state constraints are implicit functions of the decision variable (the control function). 
This section aims to clarify how one might introduce the homotopy parameter into an OCP to satisfy (A2) and (A3), thus providing one possible ``translation'' of the assumptions imposed on the NLP level to the OCP level. 

\fussy 
We consider the following parametrised discrete-time OCP, with homotopy parameter $\lambda\in[0,1]$,
\begin{numcases}{\OCP^{\lambda}:}
	\min\limits_{ u = (u_k)_{k=0}^{N-1} \subset \mbbR^m} & $\quad J(u) := \sum_{k= 0}^{N-1} \ell(x_k, u_k) + J_N(x_N)$, \nonumber \\
	\mathrm{subject\  to:}	& $\quad x_{k + 1} = f(x_k, u_k),\quad k=0,1,\dots, N - 1,$ \label{OCP_eq_1}\\
	& $\quad x_0 = x^0$,\label{OCP_eq_2}\\
	& $\quad g(\lambda, x_k) \leq \zero, \quad k=1,\dots,N,$\label{OCP_eq_3}\\
	& $\quad h(\lambda, u_k) \leq \zero, \quad k=0,\dots,N-1.$\label{OCP_eq_4}
\end{numcases} 
Here, $x_k\in\mbbR^n$ and $u_k\in\mbbR^m$ denote the state and control, respectively, at time index $k$ and $N \in \mathbb{N}$ denotes the time horizon.
The dynamics is specified by $f:\mbbR^n\times\mbbR^m\rightarrow\mbbR^n$. 
The cost functional $J:\mathbb{R}^{Nm} \rightarrow \mbbR_{\geq 0}$ consists of a running cost $\ell:\mbbR^n\times\mbbR^m\rightarrow\mbbR_{\geq 0}$ and a terminal cost $J_N:\mbbR^n\rightarrow \mbbR_{\geq 0}$.
The state and input constraints, which are perturbed by the homotopy parameter~$\lambda$, are given by $g:[0,1]\times\mbbR^n\rightarrow\mbbR^p$ and $h:[0,1]\times\mbbR^m\rightarrow\mbbR^{q}$, $p, q \in \mathbb{N}$, respectively.    
A finite-dimensional problem such as this one could be obtained from discretising the dynamics of a continuous-time optimal control problem (a direct approach), see for example \cite[Ch. 4-5]{gerdts2011optimal}. 
The reader may refer to \cite{martens2020convergence} and \cite{bonnans2017error} for results concerning convergence of the discrete-time problem's solution to the solution of the continuous-time problem, with decreasing discretisation step. 
Discrete-time OCPs are of course interesting in themselves, as many processes may be modelled by discrete-time dynamics, see for example \cite{sethioptimal2021}. 

To ease our notation, and to be consistent with previous sections, we stack the input sequence $u \subset\mbbR^m$ into one vector, $\boldu := (u_0^\top, u_1^\top, \dots, u_{N-1}^\top)^\top \in \mbbR^{r}$, $r = Nm$. Given an initial state $x^0\in\mbbR^n$, a $\boldu \in \mbbR^r$ and an index $k\in\{0,1,\dots,N\}$, we let $\Phi:\mathbb{N}_{\geq 0} \times \mbbR^r \times \mbbR^n \rightarrow \mbbR^n$ map to the state at index $k$ obtained via the dynamics, \eqref{OCP_eq_1}. That is, $\Phi(k;\boldu,x^0) =  x_k$.

We might consider more general problems that include an initial set, or where the initial state or horizon length are also decision variables. Moreover, we might want $\lambda$ to appear in other aspects of the problem, such as the cost or dynamics. Even though some of these aspects (like considering $x^0$ as a decision variable) are straightforward generalisations of the results presented here, we will not consider them as this would muddy the presentation. More difficult aspects (like $\lambda$ appearing in the dynamics) will be considered in future research.

We consider $\OCP^0$ to be an easy problem and $\OCP^{1}$ to be the difficult problem we want to solve. 
As mentioned in Section~\ref{sec:homoh_for_NLP}, our motivation for perturbing the state constraints, as in \eqref{OCP_eq_3}, comes from path planning problems through obstacle fields. 

The problem $\OCP^{\lambda}$ may of course be expressed as a nonlinear program of the form \eqref{eq:NLP}, rewritten here for convenience:
\[
\NLP^{\lambda}:\begin{cases}
	\min\limits_{\boldu\in\mathbb{R}^r} & \quad J(\boldu),\\
	\mathrm{subject\  to:} & \quad G(\lambda, \boldu) \leq \zero,
\end{cases}
\]
where $G:[0,1]\times\mbbR^r\rightarrow\mbbR^{s}$, $s = N(p + q)$, is defined follows. 
For $k = 1,2,\dots,N$ and $j = 1,2,\dots, p$,
\[
G_i(\lambda, \boldu) := g_j(\lambda,\Phi(k;\boldu,x^0)),\quad i = (j-1)N + k.
\]
For $k = 1,2,\dots,N$ and $l = 1,2,\dots, q$,
\[
G_{i}(\lambda, \boldu) := h_l(\lambda,u_{k}),\quad i = pN  + (l-1)N + k.
\]
We now impose assumptions on $\OCP^{\lambda}$.
\begin{enumerate}
	\item[\textbf{(A1$^\prime$)}] All functions appearing in $\OCP^{\lambda}$ are $C^3$ with respect to their arguments. (Those of the cost	functional, the dynamics, and the constraints.)
	\item[\textbf{(A2$^\prime$)}] There exists a $\boldu\in\mbbR^r$ such that \eqref{OCP_eq_1}-\eqref{OCP_eq_4} hold with $\lambda =1$.
	\item[\textbf{(A3$^\prime$)}] For all $\lambda\in[0,1]$, the set $\{v\in\mbbR^m : h(\lambda, v) \leq \zero\}$ is bounded.
	\item[\textbf{(A4$^\prime$)}] For $0\leq \lambda_2 \leq \lambda_1 \leq 1$ we have $G(\lambda_2, \boldu) \leq G(\lambda_1, \boldu)$, for all $\boldu\in\mbbR^r$.
\end{enumerate}
Assumption~(A2$^\prime$) says that there exists a final feasible solution to converge to, where $\lambda = 1$. 
Assumption~(A3$^\prime$) practically always holds in engineering applications. 
Assumption~(A4$^\prime$) says that as the homotopy parameter increases, the constraints become stricter. 

Our final result shows that (A2$^\prime$)-(A4$^\prime$) imply (A3), and thus that via Proposition~\ref{prop:main:NLPs} we have a convenient zero curve to track. 
\begin{proposition}\label{prop:main_result_OCP}
	Consider the problem $\OCP^{\lambda}$, $\lambda \in[0,1]$ under (A1$^\prime$) along with the homotopy $\rho$ as defined in \eqref{eq_homo}.
	For almost every $a = (\boldu^0, b^0, c^0)\in \mbbR^r \times B^0(\boldu^0) \times \mbbR^s_{>0}$, there exists a unique curve, $\Gamma_a$, emanating from the unique point $(0,\boldu^0, \boldsymbol \mu^0)$, which is $C^1$, does not intersect itself or other components of $\gamma_a$, and does not bifurcate. (See Figure~\ref{fig:zero_curves} in Section~\ref{sec:homo_for_sys_eqns}). 	
	
	Furthermore, if Assumptions (A2$^\prime$)-(A4$^\prime$) hold, along with Assumption~(A2), (A4) and (A5) of Section~\ref{sec:homoh_for_NLP}, then $\Gamma_a$ has a limit point, $(1,\boldu^\star, \boldsymbol \mu^\star)$, where the pair $(\boldu^\star, \boldsymbol \mu^\star)$ is a local solution to $\NLP^1$, and thus to $\OCP^1$. 
	If, additionally, (A6) holds, then $\Gamma_a$ has finite arc length. 
\end{proposition}
\begin{proof}
	The first part of the proposition follows directly from Proposition~\ref{prop:main:NLPs}. 
	We now argue that (A2$^\prime$)-(A4$^\prime$) imply (A3). 
	
	Assumption (A4$^\prime$) implies that, for any $b \in\mbbR_{> 0}$,
	\[
	G(\lambda_2, \boldu) \leq G(\lambda_1, \boldu) \leq (1 - \lambda_1)b \leq (1 - \lambda_2)b,\quad 0\leq \lambda_2 \leq \lambda_1 \leq 1.
	\] 
	Moreover, recall that $S(\lambda, b) = \{\boldu : G(\lambda, \boldu) \leq (1 - \lambda) b\}$, for $b \in\mbbR_{>0}$ and $\lambda \in[0,1]$ (see \eqref{def:S}). 
	Therefore, with (A2$^\prime$) we get,
	\[
	S(1,b) = \{\boldu : G(1,\boldu)\leq \zero\} \neq \emptyset.
	\]
	Therefore, for any $b^0\in B^0(\boldu^0)$, where $\boldu^0\in \mbbR^r$,
	\[
	\emptyset \neq S(1, b^0) \subseteq S(\lambda_1, b^0) \subseteq S(\lambda_2, b^0) \subseteq S(0, b^0).
	\]
	Under (A3$^\prime$) as well, we get that $\cup_{\lambda\in[0,1]} S(\lambda, b^0)$ is bounded. 
	Therefore, (A3) holds, and we may deduce that if (A2), (A4) and (A5) also hold, then $\Gamma_a$ has the stated properties via Proposition~\ref{prop:main:NLPs}. 
	As before, the final statement involving (A6) follows from Lemma~\ref{lem:homo_is_nice}. 
\end{proof}

\subsection*{Comments on the main result, Proposition~\ref{prop:main_result_OCP}}

This section aims to clarify the results of Proposition~\ref{prop:main_result_OCP}, and to show how it may be used to solve difficult optimal path planning problems, in particular. 

First, decide how the homotopy parameter should be introduced such that Assumptions (A2) and (A4$^\prime$) hold. 
Thus, one needs to identify the constraints that are nonconvex with respect to $\boldu$, and introduce $\lambda$ such that these constraints disappear when $\lambda = 0$ (A2), and become stricter with increasing $\lambda$ (A4$^\prime$). 

Next, consider \emph{any} initial vector, $\boldu^0\in\mbbR^r$, along with any $b^0\in B^0(\boldu^0)$ and any $c^0\in\mbbR_{> 0}$.  
Note that $\boldu^0$ does not even need to be a feasible solution to the relaxed problem, where $\lambda = 0$. (Recall that the homotopy $\rho$ is \emph{globally} convergent). 
However, if $\boldu^0$ satisfies $G(0,\boldu^0) \leq \zero$, which is an easy convex problem by (A2), then $B^0(\boldu^0) = \mbbR^s_{> 0}$. 

Next, with $a = (\boldu^0, b^0, c^0)$, find the unique $\boldsymbol \mu^0\in\mbbR^s$ for which $\rho_a(0,\boldu^0,\boldsymbol \mu^0) = \zero$. 
Looking at the definition of $K$, see \eqref{eq:K}, this amounts to solving for each $\mu_i$, $i\in[1:s]$, in the equation,
\[
	\mu_i^3 - | C_{i1} - \mu_i |^3 + C_{i2} = 0,
\]
where $C_{i1},C_{i2}\in\mbbR$ are constants. 
The final step is to track $\Gamma_a$ numerically, starting from $(0,\boldu^0,\boldsymbol \mu^0)$, for which there are many techniques, see for example \cite{zangwil1981pathways, allgower2003introduction}. 
If the remaining assumptions of Proposition~\ref{prop:main_result_OCP} hold and the curve tracking is successful one will arrive at a local solution to $\OCP^1$, namely $(\boldu^\star, \boldsymbol \mu^\star)$.

For path planning problems it is typical that the obstacles are nonconvex in the state, $x_k$, but that the target set and control constraints are convex in $x_k$ and $u_k$, respectively. 
Thus, a sensible approach for the first step is to introduce $\lambda$ such that it shrinks the obstacles to points when $\lambda = 0$, and grows them, to their original sizes, as $\lambda$ increases.
For linear dynamics, that is, $f(x_k,u_k) = Ax_k + Bu_k$, with $A_k\in\mbbR^{n\times n}$, $B_k\in\mbbR^{n\times m}$, $A_k \neq 0$ and $B_k\neq 0$, a constraint that is convex (resp. nonconvex) in the state remains convex (resp. nonconvex) when expressed as a function of the input sequence, $\boldu$. 
Thus, identifying $\mathbb{I}_{nc}$ in (A2) is easy in this special case.

\section{Numerical examples}\label{sec:numerics}

Although our results are applicable to optimal control problems in general, we focus on optimal path planning problems, which were the main motivation behind this paper's research. 
We first provide details of a curve tracking algorithm and then solve a number of problems that highlight the strengths of the approach. 

\subsection{Curve tracking algorithm}

\begin{figure}[htb]
	\centering
	\includegraphics[width=0.6\linewidth]{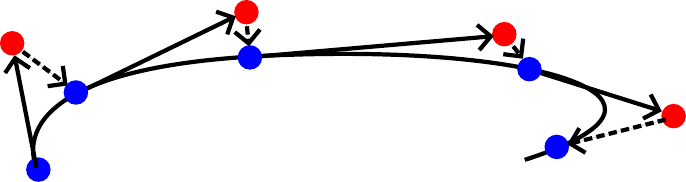}
	\caption{Intuitive idea of a predictor-corrector algorithm. 
		The predictor generates points (red dots), following the tangents to the curve at the current iterate. The corrector determines the next point on the curve (blue dots) with a root finding algorithm (like Newton's method) initiating from the predicted point.}
	\label{fig:curve_tracking_idea}
\end{figure}

We implemented a path tracking algorithm in Matlab, see \cite{watson1997algorithm} for a Fortran package. 
The important steps of the algorithm are to obtain a vector tangent to the zero curve pointing in a consistent direction (the \textbf{getTangent} function); to obtain a point along this tangent direction (the ``predict'' step); and to then find a point on the zero curve close to the predicted step (the ``correct'' step, performed by the \textbf{corrector} function). 
See Figure~\ref{fig:curve_tracking_idea} for some intuition. 

Following \cite[Ch.4]{allgower2003introduction}, Lines~20~-~21 of the \textbf{getTangent} function obtains a direction that is tangent to the zero curve via the Q-R decomposition of $\nabla \rho_a(w^i)^\top$. 
In more detail, let
\[
\nabla \rho_a(w^i)^\top = Q \left(\begin{array}{c} R\\ \mathbf{0}^\top \end{array} \right),
\]
where $Q$ is orthogonal and $R$ is upper-triangular. 
Then $t = Q_{[\bigcdot,\mathrm{end}]}$ (the last column of $Q$) is a vector $t\in\mbbR^{r+s+1}$ satisfying $\nabla \rho_a(w^i)^\top t = 0$, and $\Vert t \Vert = 1$. 
If $\det\begin{pmatrix}\nabla \rho_a(w^i)\\ t^\top \end{pmatrix}~>~0$ (resp. $\det\begin{pmatrix}\nabla \rho_a(w^i)\\ t^\top \end{pmatrix}~<~0$) the vector $t$ points in the direction of increasing (resp. decreasing) arc length. 
On the first iteration of the \textbf{while} loop (i.e, for $i=0$) it might not be clear in which direction the zero curve should be followed. 
The \textbf{if} statement in Lines~22~-~30 checks this and makes sure the direction remains consistent throughout iterations of $i$.
The functions \textbf{corrector} and \textbf{finalCorrector} implement Newton's method. 
The notations $\rho_a(\bigcdot)$ and $\nabla \rho_a(\bigcdot)$ are meant to indicate that functions are passed to these correctors. 
Algorithm~\ref{alg:1} is of ``Euler-Newton'' type, and is guaranteed to trace out the homotopy's zero curve with a sufficiently small $h>0$, see \cite[Thm.~5.2.1]{allgower2003introduction}. 
\begin{algorithm}
	\caption{}
	\label{alg:1}
	\begin{algorithmic}[1]
		\State{Choose a $\boldu^0 \in \mathbb{R}^r$.}
		\State{Choose a $b^0 \in \mathbb{R}_{> 0}^s$ such that $G(0,\boldu^0) \leq b^0$, and $c^0 \in \mathbb{R}^s_{>0}$.}
		\State{Let $a \gets (\boldu^0, b^0, c^0)$.}
		\State Get symbolic expressions for $\rho_a(\bigcdot)$ and $\nabla \rho_a(\bigcdot)$.
		\State{Find the unique $\boldsymbol \mu^0 \in \mathbb{R}^s_{>0}$ for which $\rho_a(0,\boldu^0,\boldsymbol\mu^0) = \mathbf{0}$.}
		\State{Let $i \gets 0$, $\boldu^i \gets \boldu^0$, $\boldsymbol \mu^i \gets \boldsymbol\mu^0$, $\lambda^i \gets 0$.}
		\State{Let $w^i \gets (\lambda^i, \boldu^i, \boldsymbol \mu^i)$.}
		\State{Choose arc length step size $h > 0$.}
		\State{Choose a tolerance $\mathrm{tol}_{\mathrm{main}} > 0$ for the termination criterion.}
		\While{$| \lambda^i - 1 | < \mathrm{tol}_{\mathrm{main}}$}
		\State{$t \gets$ \textsc{getTangent}$(\nabla\rho_a(w^i))$ \texttt{// get vector tangent to the zero curve.}}
		\State{$w^{i+1} \gets w^i + h \cdot t$ \texttt{// Predict.}}
		\State{$w^{i+1}\gets$ \textsc{corrector}$(w^{i+1}, \rho_a(\bigcdot), \nabla \rho_a(\bigcdot))$ \texttt{// Correct.}}
		\State{Let $i \gets i + 1$}
		\EndWhile
		\State{\texttt{// Final corrector step}}
		\State Let $(\boldu^\mathrm{sol}, \boldsymbol \mu^\mathrm{sol})\gets$\textsc{finalCorrector}$(w^i,\rho_a(\bigcdot))$
		\State
		\Function{getTangent}{$\nabla\rho_a(w^i)$}
		\State{Let $Q \left(\begin{array}{c} R\\ \mathbf{0}^\top \end{array} \right) \gets \textsc{getQRDecomposition}(\nabla\rho_a(w^i)^\top)$}
		\State{Let $t \gets Q_{[\bigcdot,\mathrm{end}]}$} \texttt{// Let $t$ be the last row of $Q$.}
		\If{$\lambda$ initially increases with increasing arc length}
			\If{$\det(Q)\det(R) < 0$} 
				\State{Let $t\gets -t$}
			\EndIf
		\Else
			\If{$\det(Q)\det(R) >0$} 
				\State{Let $t\gets -t$}
			\EndIf
		\EndIf
		\State \Return $t$
		\EndFunction
		\State
		\Function{corrector}{$w^{i+1}, \rho_a(\bigcdot), \nabla \rho_a(\bigcdot)$}
		\State{Choose $\mathrm{maxIterations} > 0$ and $\mathrm{tol}_{\mathrm{Newton}} > 0$} 
		\State{Let $w^{\mathrm{loop}}\gets w^{i+1}$}
		\State{Let $j\gets 0$}
		\While{$(\Vert \rho_a(w^{\mathrm{loop}}) \Vert > \mathrm{tol}_{\mathrm{Newton}}) \mathrm{\,\,\bf{or}\,\,} (j < \mathrm{maxIterations})$ }
		\State{Let $w^{\mathrm{loop}} \gets w^{\mathrm{loop}} - \nabla \rho_a^{\dagger}(w^{\mathrm{loop}}) \rho_a(w^{\mathrm{loop}})$}
		\EndWhile
		\State \Return $w^{\mathrm{loop}}$
		\EndFunction
		\State
		\Function{finalCorrector}{$w^{i}, \rho_a(\bigcdot)$}
		\State{Choose $\mathrm{maxIterations} > 0$ and $\mathrm{tol}_{\mathrm{Newton}} > 0$} 
		\State Let $\rho_a^{\mathrm{final}}(\bigcdot, \bigcdot) \gets \rho_a(1,\bigcdot, \bigcdot)$
		\State{Let $(\boldu^{\mathrm{loop}},\boldsymbol \mu^{\mathrm{loop}})\gets (w^{i}_{[1:r]}, w^i_{[r+1:r+s]})$}
		\State{Let $j\gets 0$}
		\While{$(\Vert \rho_a^{\mathrm{final}}(\boldu^{\mathrm{loop}},\boldsymbol\mu^{\mathrm{loop}})) \Vert > \mathrm{tol}_{\mathrm{Newton}}) \mathrm{\,\,\bf{or}\,\,} (j < \mathrm{maxIterations})$ }
		\State{$(\boldu^{\mathrm{loop}},\boldsymbol\mu^{\mathrm{loop}}) \gets (\boldu^{\mathrm{loop}},\boldsymbol\mu^{\mathrm{loop}}) - \nabla_{\boldu,\boldsymbol\mu} \rho_a^{\mathrm{final}}(\boldu^{\mathrm{loop}},\boldsymbol\mu^{\mathrm{loop}})^{-1} \rho_a^{\mathrm{final}}(\boldu^{\mathrm{loop}},\boldsymbol\mu^{\mathrm{loop}})$}
		\EndWhile
		\State \Return $(\boldu^{\mathrm{loop}},\boldsymbol\mu^{\mathrm{loop}})$
		\EndFunction
	\end{algorithmic}
\end{algorithm}

\subsection{Linear system}

Consider the following problem, where $x_k,u_k\in\mbbR^2$,
\begin{align}
	\min\limits_{(u_k)_{k=0}^{N-1}\subset\mbbR^2} & \quad \sum_{k= 0}^{N-1} \frac 12 \Vert u_k \Vert^2 + \frac 12 \beta \| x_N - x^{\mathrm{Target}} \|^2 , \label{eq:linear_example_1} \\
	\mathrm{subject\  to:}	& \quad x_{k+1}  = x_k + u_k, \quad k=0,1,\dots,N-1\label{eq:linear_example_2} \\
	& \quad x_0 = x^0, \label{eq:linear_example_3} \\
	& \quad-\Vert x_k - m_1 \Vert^2 + \lambda r  \leq 0,\,\, k=1,\dots,N,	\label{eq:linear_example_4} \\
	& \quad-\Vert x_k - m_2 \Vert^2 + \lambda r \leq 0 ,\,\, k=1,\dots,N,	\label{eq:linear_example_5} \\
	&\quad \Vert u_k \Vert^2  \leq 1,\,\, k=0,1,\dots,N-1. \label{eq:linear_example_6}
\end{align}
The problem is to find a path that starts from $x^0\in\mbbR^2$, ends close to $x^{\mathrm{Target}}\in\mbbR^2$, avoids two circles centred at $m_1\in\mbbR^2$ and $m_2\in\mbbR^2$ of radius $r \geq  0$, while minimising a compromise between minimal control effort and distance to the target. 
The constant $\beta >0$ may be used to weight the importance of these two objectives. 

Referring to Proposition~\ref{prop:main_result_OCP}, all functions are smooth, thus (A1$^\prime$) holds. 
Because the dynamics is simple it is easy to specify a path that is feasible for the desired problem where $\lambda = 1$, so (A2$^\prime$) is also satisfied. 
The control is bounded, so (A3$^\prime$) holds, and we have introduced $\lambda$ such that with $\lambda = 0$ the obstacles are simply points and grow with increasing $\lambda$. 
Thus, (A2) and (A4$^\prime$) hold. 
We can deduce from Proposition~\ref{prop:main_result_OCP} that there exists a convenient zero curve, $\Gamma_a$, for almost every $a$. 

\begin{figure}[h]
	\centering
	\includegraphics[scale=0.6]{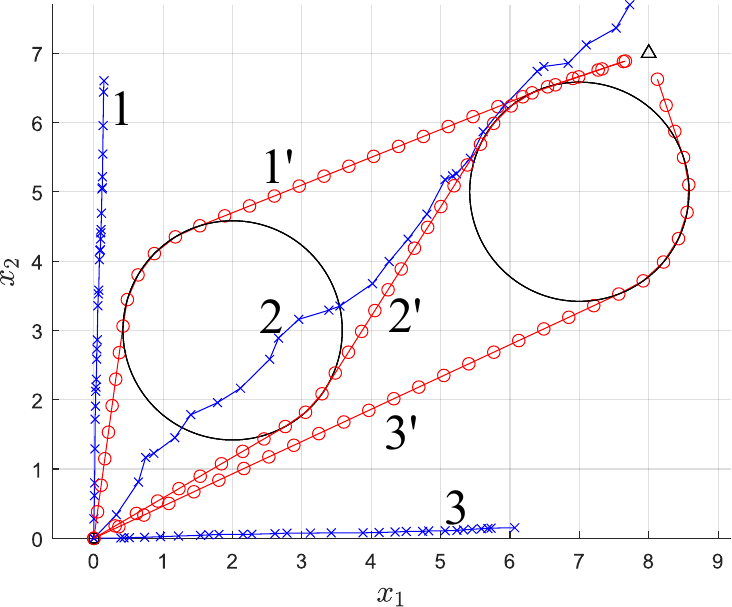}
	\caption{Results of tracking the zero curve of the parametrised homotopy, $\rho_a$, see \eqref{def:para_homo}, using Algorithm~\ref{alg:1}, for the linear problem, \eqref{eq:linear_example_1}-\eqref{eq:linear_example_6}. Initial paths, associated with guesses $\boldu^0$, are indicated by the blue crosses (curves 1,2 and 3). The local solutions the algorithm converges to are indicated by red circles (curves 1', 2', 3').}
	\label{fig:Linear_2_obs_no_overlap_a}
\end{figure}

\begin{figure}[h]
	\centering
	\includegraphics[scale=0.6]{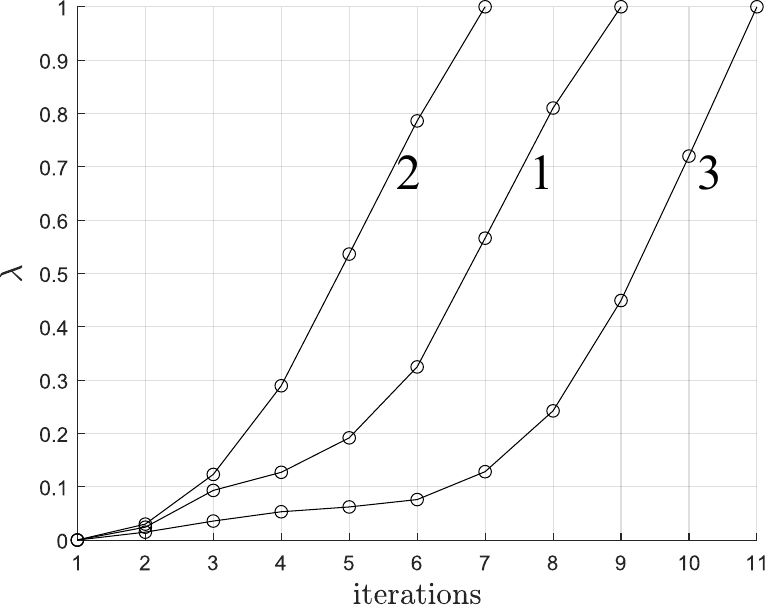}
	\caption{Homotopy parameter over iterations for the problem in Figure~\ref{fig:Linear_2_obs_no_overlap_a}. Numbers indicate which initial guess, $\boldu^0$, the sequence is associated with.}
	\label{fig:Linear_2_obs_no_overlap_b}
\end{figure}

We consider the problem with $m_1 = (2, 3)^\top$, $m_2 = (7, 5)^\top$, $r=2$, $x^0 = (0,0)^\top$, $x^{\mathrm{Target}} = (8,7)^\top$, $N = 30$, and $\beta=1$. 
Thus, referring to the beginning of Section~\ref{sec:homo_for_OCPs}, we have $r = Nm = 60$ and $s = N(p+q) = 90$, so we need to track the zero curve, $\Gamma_a$, which lives in $r + s =150$ dimensions. 
We invoke Algorithm~\ref{alg:1} with step size $h = 0.5$. 
Figure~\ref{fig:Linear_2_obs_no_overlap_a} shows various state paths, indicated by blue crosses, each resulting from a different initial guess, $\boldu^0\in\mbbR^r$. 
The algorithm successfully tracks the zero curve to $\lambda = 1$ for each of these $\boldu^0$'s, converging to locally optimal solutions, indicated with red circles. 
Note how the local solution obtained depends on the initial guess. 
Figure~\ref{fig:Linear_2_obs_no_overlap_b} shows how $\lambda$ changes as $i$ increases in Algorithm~\ref{alg:1}, clearly demonstrating that a better initial guess at the solution results in a shorter path to track. 
The calculation of the symbolic Jacobian $\nabla \rho_a(\bigcdot)$ (which we did with CasADi, \cite{Andersson2019}) takes up the majority of the algorithm's computational effort. 
On a mid-range PC the average time to calculate $\nabla \rho_a(\bigcdot)$ over 100 random $\boldu^0$'s was $\sim 6.25$ seconds, whereas the average time to track the zero curve to the solution was $\sim0.3$ seconds. 

\begin{figure}[h]
	\centering
	\includegraphics[scale=0.6]{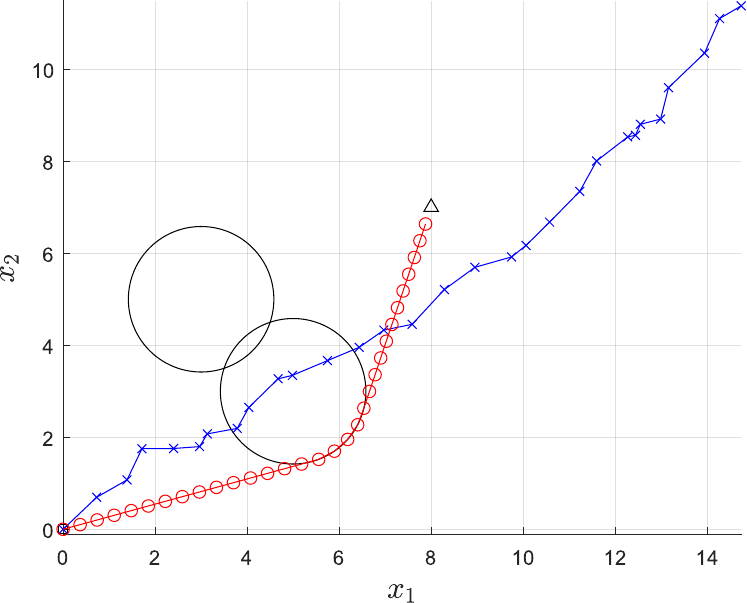}
	\caption{Result of tracking the zero curve for the linear problem, \eqref{eq:linear_example_1}-\eqref{eq:linear_example_6} with the obstacles moved so they overlap. 
	Blue crosses: initial guess, red circles: solution.}
	\label{fig:Linear_overlapping_converges}
\end{figure}

\begin{figure}[!h]
	\centering
	\begin{subfigure}[t]{0.45\linewidth}
		\centering
		\includegraphics[scale=0.44]{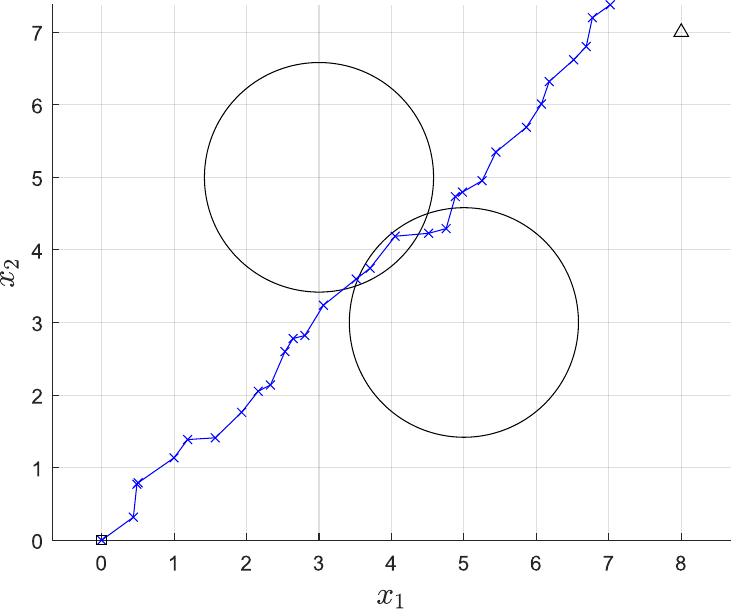}
		\caption{The algorithm does not converge with this initial $\boldu^0$.}
		\label{fig:Linear_2_obs_overlap_a}
	\end{subfigure}\hfill
	\begin{subfigure}[t]{0.45\linewidth}
		\centering
		\includegraphics[scale=0.44]{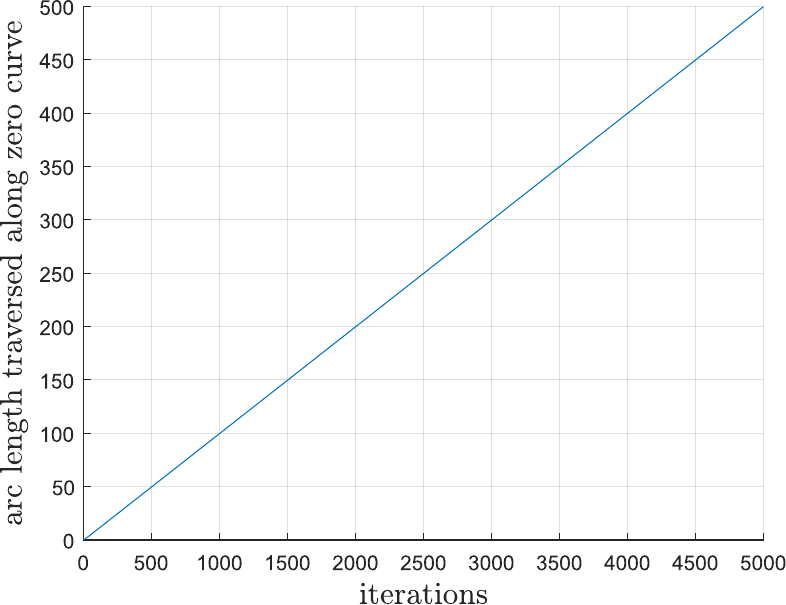}
		\caption[b]{The arc length strictly increases suggesting that $\Gamma_a$ has infinite arc length.}
		\label{fig:Linear_2_obs_overlap_b}
	\end{subfigure}
	\caption{Tracking the parametrised zero curve, $\Gamma_a$, using Algorithm~\ref{alg:1}, for the linear problem, \eqref{eq:linear_example_1}-\eqref{eq:linear_example_6}, with control bound $0.5$, and $\beta=20$ in the cost function.}
	\label{fig:Linear_2_obs_overlap_1}
\end{figure}

We now adapt the linear example. 
We change the obstacle centres and radii so that they overlap with $\lambda = 1$, as shown in Figure~\ref{fig:Linear_overlapping_converges}. 
While the algorithm often converges, even for initial guesses passing between the two obstacles when they are shrunk, there are instances where it fails. 
Figure~\ref{fig:Linear_2_obs_overlap_1} shows such an example, where we decrease the bound on the norm of the control, as in \eqref{eq:linear_example_6}, to $0.5$ and select $\beta = 20$ in the cost function. 
The homotopy parameter, $\lambda$, remains strictly less that 1 over thousands of iterations. Figure~\ref{fig:Linear_2_obs_overlap_b} shows the traversed arc length along the zero curve against increasing iterations in the algorithm, with a smaller step length of $h=0.1$. 
For large $i$ the Jacobian $\nabla \rho_a(1, \boldu^{i} , \boldsymbol  \mu^{i})$ is close to singular, and the final Newton corrector invoked in Line~17 of the algorithm always fails if invoked at some large $i$. 
These facts would suggest that (A6) is violated for this choice of $\boldu^0$, and that we are tracking a zero curve with infinite arc length, leading towards a degenerate KKT point. 
This example emphasises that one may still need to choose $a$ carefully, even though there exists a zero curve for almost every $a$. 

\subsection{Nonlinear example: Dubins' vehicle}

\begin{figure}[h]
	\centering
	\includegraphics[scale=0.6]{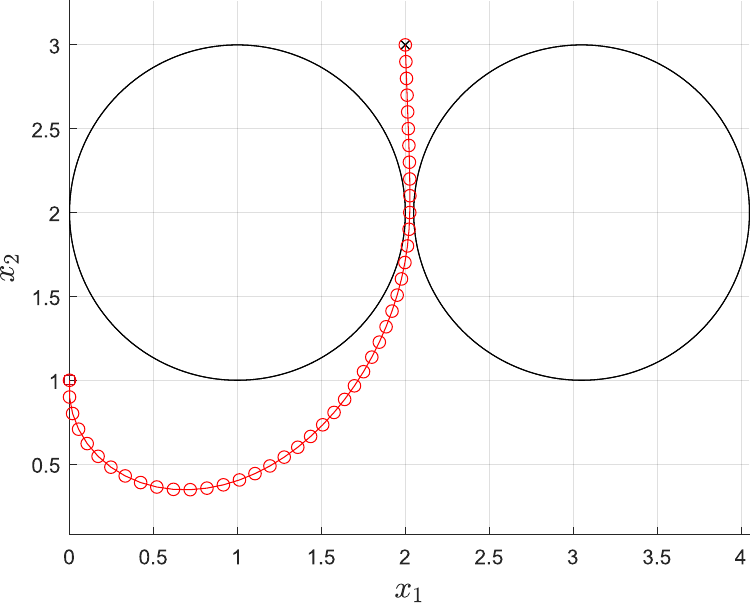}
	\caption{Optimal solution obtained for the OCP, \eqref{eq:dubins1}-\eqref{eq:dubins8}, with Algorithm~\ref{alg:1}.}
	\label{fig:Dubins_homo}
\end{figure}

We now demonstrate that the approach is just as applicable to nonlinear systems. 
Consider the continuous-time optimal control problem involving the Dubins vehicle, \cite{dubins1957curves},
\begin{align}
	\min\limits_{u\in\mathcal{L}_T} & \quad \Vert \x(T) - \x^{\mathrm{Target}} \Vert^2\, \label{eq:dubins1} \\
	\mathrm{subject\  to:}	& \quad \dot x_1(t) = v\cos \theta(t), \label{eq:dubins2}\\
	& \quad \dot x_2(t) = v\sin \theta(t), \label{eq:dubins3}\\
	& \quad \dot \theta(t) = u(t), \label{eq:dubins4}\\
	& \quad (x_1(0), x_2(0), \theta(0))^\top = (\x^0, \theta^0)^\top,\label{eq:dubins5} \\
	& \quad-\Vert \x(t)- m_1 \Vert^2+\lambda r_1   \leq 0, \label{eq:dubins6}\\
	& \quad-\Vert \x(t) - m_2 \Vert^2 + \lambda r_2 \leq 0,\label{eq:dubins7}\\
	&\quad | u(t) | \leq u^{\max},\label{eq:dubins8}
\end{align}
$t\geq 0$. 
Here, $\x = (x_1,x_2)\in\mbbR\times\mbbR$ is the car's position $\theta\in\mathbb{S}$ is its orientation and $v\in\mbbR_{>0}$ is its constant speed. 
Its turning rate is determined by the control, $u\in\mbbR$. 
It is easily verified that with a constant $u$  and $v$ the car traces out circles in the $x_1$--$x_2$ space, and that its minimal turning radius is $\frac{v}{u^{\max}}$. 
By $\mathcal{L}_T$ we denote Lebesgue-measurable functions mapping $[0,T]$ to $\mbbR$. 
We choose the centres, $m_i$, and radii, $r_i$, to create the obstacle field in Figure~\ref{fig:Dubins_homo}, and select the initial and target states to be $(x_1^0,x_2^0,\theta^0) = (0,1,-\pi/2)^\top$ and $(x_1^{\mathrm{target}},x_2^{\mathrm{target}},\theta^{\mathrm{target}}) = (2,3,\mathrm{free})^\top$, respectively.  
Concerning the assumptions, as in the linear example all functions are $C^3$; the nonconvex constraints disappear with $\lambda = 0$; the feasible space becomes stricter with increasing $\lambda$; and the control is bounded. 
Because the dynamics is not simple anymore, we cannot easily verify that a feasible solution exists. 

Nevertheless, we select $v = 1$ and $u^{\max} = 6$ in the model and run Algorithm~\ref{alg:1} with the parameters $h=0.25$, $N=44$, discretising the dynamics via the Euler scheme with a step size of $\delta t =0.1s$.  
We consistently obtain the solution indicated in Figure~\ref{fig:Dubins_homo}, regardless of the random intial guess, $\boldu^0$. 
On the same PC as before, the time taken to compute the symbolic Jacobian, $\nabla \rho_a(\bigcdot)$, was about 22~seconds, whereas it took about 1.6~seconds to track the zero curve to the solution, over 30 predictor-corrector steps.  

\subsection{Comparison with sampling-based motion planners}

We now briefly compare the convergent homotopy approach with sampling-based motion planners, see \cite[Ch.~5]{lavalle2006planning}. 
These algorithms take the following general steps: 
\begin{enumerate}
	\item Consider an initial graph, $\mathcal{G}(V,E)$, containing the initial and target states, $\x^{0}$ and $\x^{\mathrm{target}}$ as vertices.
	\item Select a vertex $\x^{\mathrm{current}}\in\mathcal{G}(V,E)$. 
	\item Select a state from the free space, $\x^{\mathrm{new}}\in\mathcal{C}_{\mathrm{free}}$, (i.e., not inside an obstacle) and attempt to construct a collision-free path between $\x^{\mathrm{current}}$ and $\x^{\mathrm{new}}$, i.e., a curve $\tau:[0,1]\rightarrow \mathcal{C}_{\mathrm{free}}$ such that $\tau(0)=\x^{\mathrm{current}}$ and $\tau(1)=\x^{\mathrm{new}}$. If this is not possible, go back to Step 2).
	\item Add $\tau$ as an edge of $\mathcal{G}$.
	\item Check if there exists a path from $\x^{0}$ to $\x^{\mathrm{target}}$. If not, go to Step 2).
\end{enumerate}
Our homotopy approach has some advantages over sampling-based planners. 
First, it does not suffer from the ``curse of (state) dimensionality''. 
In some sampling-based planners, such as A*, one must create a grid of the state space, with the number of points required in this sampling growing exponentially with the state dimension, see \cite{karaman2011sampling}. 
On the other hand, in Algorithm~\ref{alg:1} the computational effort increases independent of the state dimension, as the problem is to track a curve in $N(m+p+q)$ dimensions. 

\begin{figure}[h]
	\centering
	\includegraphics[scale=0.6]{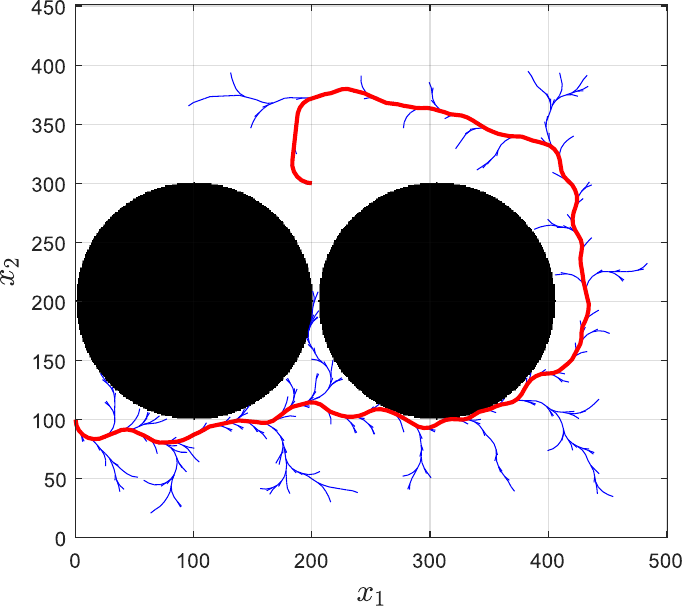}
	\caption{The thick red curve is the first feasible solution found to OCP \eqref{eq:dubins1}-\eqref{eq:dubins8} via Matlab's $\mathsf{plannerRRTStar}$. 
	The thin blue curves are areas searched by the algorithm.}
	\label{fig:Dubins_RRTStar}
\end{figure}

Second, sampling-based algorithms may struggle to find paths through complicated obstacle fields, especially if there are narrow passages. 
Figure~\ref{fig:Dubins_RRTStar} shows the first feasible path found with the Matlab solver $\mathsf{plannerRRTStar}$, with appropriate parameters so that the curves obtained in the $x_1$--$x_2$ space correspond to those of a Dubins car, travelling at unit speed with a minimum turning radius of $\frac{1}{6}$, as in the problem \eqref{eq:dubins1}-\eqref{eq:dubins8}. 
Matlab's $\mathsf{plannerRRTStar}$ implements the popular RRT* algorithm, introduced in \cite{karaman2011sampling}, which is \emph{asymptotically optimal}, meaning that it will find the shortest path with infinite running time. 
Although the time taken to find the first feasible solution is very fast (roughly 1~second on the same computer) the path obtained depends on the algorithm's random searching pattern, which we found often misses the shorter path passing between the two obstacles. 
This solution rarely improves, even if the algorithm is allowed to run for many thousands of extra iterations. 
Moreover, the paths found via sampling-based algorithms are often jagged and correspond to control signals that are erratic, whereas those found with the homotopy approach are often of higher quality. 

\begin{figure}[h]
	\centering
	\includegraphics[scale=0.6]{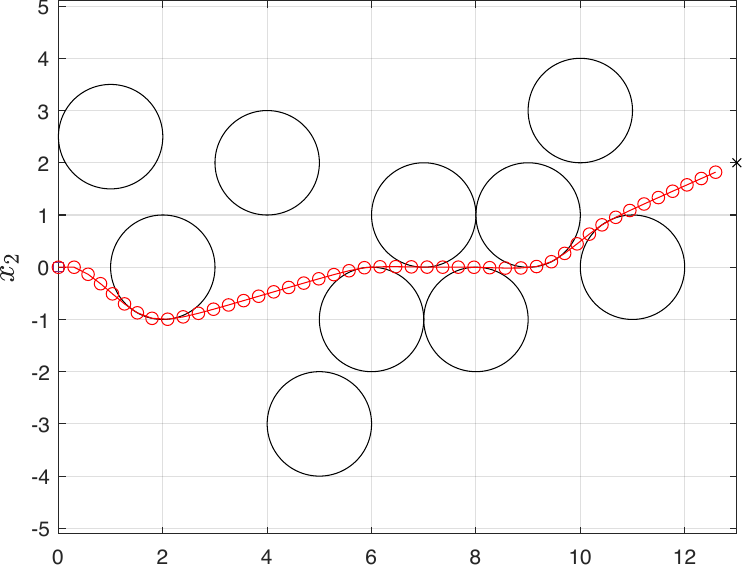}
	\caption{Solution found for the Dubins vehicle, \eqref{eq:dubins2}-\eqref{eq:dubins4}, from the initial point, $(0,0,0)^\top$, to the target, $(13,2,\mathrm{free})$, through an obstacle field.}
	\label{fig:Dubins_many_obs_homo}
\end{figure}

\begin{figure}[h]
	\centering
	\includegraphics[scale=0.6]{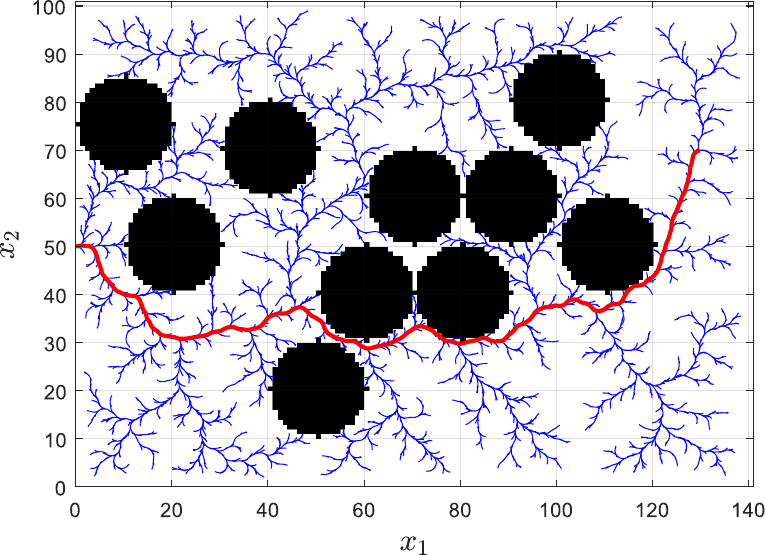}
	\caption{First feasible solution (thick red curve) found to the same problem as in Figure~\ref{fig:Dubins_many_obs_homo} with Matlab's $\mathsf{plannerRRTStar}$.}
	\label{fig:Dubins_RRTSTAR_many_obs}
\end{figure}

As a final example, Figure~\ref{fig:Dubins_many_obs_homo} shows the solution found to a problem with many obstacles via the homotopy method. 
On the same PC, the time taken to compute the Jacobian, $\nabla \rho_a$, was about 53~seconds and it took about 23~seconds to converge to the solution, in 128 predictor-corrector steps. 
Figure~\ref{fig:Dubins_RRTSTAR_many_obs} shows the first feasible solution found with $\mathsf{plannerRRTStar}$, which took 23 seconds. 
The increase in computation time for our approach results in the long horizon and large number of obstacles: the algorithm needs to track a curve in $N(m + p + q) = 46(1+10+1) = 552$ dimensions. 


From our timings, a valid point to raise is that the computational complexity of the homotopy approach increases rapidly with a large number of obstacles, whereas an algorithm like RRT* performs well in such cases because it does not explicitly construct the obstacle field, but merely probes it through sampling. 
However, note that the vast majority of computational effort goes towards pre-computing the symbolic Jacobian matrix, $\nabla \rho_a(\bigcdot)$, whereas the time taken to track the zero curve is often much faster. 
These facts might be exploited to have a reliable fast solver in model predictive control (MPC), see \cite{boccia2014stability, grune2017nonlinear}. 
In this approach a finite-horizon OCP is iteratively solved with a receding horizon to obtain a closed-loop controller. 
If the initial state, along with the control, is included as a decision variable in the finite-horizon OCP, then, with a particular fixed $a$, the symbolic Jacobian, $\nabla\rho_a$, can be precomputed and successively reused to solve the OCP.  
To use MPC online in a robotic application, for example, one would utilise a prediction horizon that is much shorter that the 4.6 seconds used in the problem shown in Figure~\ref{fig:Dubins_many_obs_homo}, thus the computation time might be cut down considerably (recall that tracking the zero curve for the problem in Figure~\ref{fig:Dubins_RRTStar}, which lives in 220 dimensions, took just 1.6 seconds). 
Finally, results from the paper \cite{flasskamp2019symmetry}, which deals with how motion primitives and MPC can be combined, might be helpful to construct a good initial guess, $\boldu^0$.

\section{Conclusion}\label{sec:conclusion}

We proposed novel sufficient conditions on discrete-time optimal control problems under which a homotopy can be rendered globally convergent with probability one. 
To this end, we leveraged the key results presented in~\cite{watson2001theory,chow1978finding}, for NLPs.  
We demonstrated the applicability of our results by solving non-convex optimal path planning problems by tracking the homotopy's zero curve via a numerical scheme. 
We also pointed out the possible use of the approach as an online solver for model predictive control, if the Jacobian $\nabla\rho_a$ is pre-computed offline. 

In future work, one could study the (measure-zero) set of initial parameters $a = (\boldu^0,b^0,c^0)$ 
for which no convenient zero path exist. 
Then one may investigate the physical interpretation of these singularities, e.g., in the context of robotic path planning. 
Also being able to include the homotopy parameter in the dynamics and/or cost function is an important generalisation that should be considered.  

We think that our initial numerical experiments show some promise for the approach. 
An in-depth numerical study could be conducted comparing the homotopy approach of this paper (implemented with faster adaptive step-size algorithms) with the plethora of other robotic motion planning algorithms, comparing execution time and quality of solutions found. 
Moreover, the approach could also be compared with other state-of-the-art solvers of nonlinear and nonconvex optimisation problems.
The claim we made that the approach could be attractive for online MPC of robots could also be rigorously investigated. 

Finally, although a lot of research has been done on effective path-tracking algorithms, see for example \cite{zangwil1981pathways,allgower2003introduction}, Matlab and Python code implementing these techniques are still scarce
such that toolboxes in these languages might be desirable.

\section*{Acknowledgments}
We thank Manuel Schaller for excellent comments on how to improve the manuscript.

\newpage
\bibliographystyle{IEEEtran}
\bibliography{homotopy_methods_in_OC}

\end{document}